\documentclass[11pt,a4paper,reqno]{amsart}

\usepackage[utf8]{inputenc}
\usepackage[T1]{fontenc}
\usepackage{xspace,
amsthm,
amsmath,
amssymb,
bbm,
tikz,
graphicx,
subfig,
xcolor,
keyval}
\DeclareGraphicsExtensions{.png, .pdf}

\usepackage[hidelinks]{hyperref}

\usepackage{geometry}
\tikzstyle{nodino}=[circle,draw,fill,inner sep=0pt,minimum size=0.5mm]
\tikzstyle{infinito}=[circle,inner sep=0pt,minimum size=0mm]
\tikzstyle{nodo}=[circle,draw,fill,inner sep=0pt, minimum size=0.5*width("k")]
\tikzstyle{nodo_vuoto}=[circle,draw,inner sep=0pt, minimum size=0.5*width("k")]
\usetikzlibrary{graphs}
\tikzset{every loop/.style={min distance=10mm,in=300,out=240,looseness=10}}
\tikzset{place/.style={circle,thick,draw=blue!75,fill=blue!20,minimum size=6mm}}
\tikzset{place2/.style={circle,thick,draw=red!75,fill=red!20,minimum size=6mm}}

\newcommand{\f}[2]{\frac{#1}{#2}}
\newcommand{\tf}[2]{\tfrac{#1}{#2}}
\newcommand{\E}{\mathrm{E}}
\newcommand{\V}{\mathrm{V}}
\newcommand{\G}{\mathcal{G}}
\newcommand{\F}{\mathcal{F}}
\newcommand{\M}{\mathcal{M}}
\newcommand{\D}{\mathcal{D}}
\newcommand{\K}{\mathcal{K}}
\renewcommand{\L}{\mathcal{L}}
\newcommand{\N}{\mathbb{N}}
\newcommand{\R}{\mathbb{R}}
\newcommand{\C}{\mathbb{C}}
\newcommand{\Q}{\mathcal{Q}_{\Dg}}
\newcommand{\lap}{\Delta_\G}
\newcommand{\Dg}{\D_\G}
\newcommand{\dx}{\,dx}
\newcommand{\Ci}{\mathcal{C}(\G,\infty)}
\newcommand{\Cp}{\mathcal{C}(\G,p)}
\newcommand{\Cs}{\mathcal{C}(\G,6)}
\newcommand{\Ck}{\mathcal{C}(\K)}
\renewcommand{\H}{\mathcal{H}(\G)}
\newcommand{\Hm}{\mathcal{H}_\mu(\G)}
\renewcommand{\d}{\mathrm{dom}}
\renewcommand{\v}{\mathrm{v}}
\newcommand{\ep}{\varepsilon}

\theoremstyle{plain} 
\newtheorem{thm}{Theorem}[section]

\theoremstyle{definition}

\theoremstyle{definition}

\theoremstyle{remark} 
\newtheorem{rem}{Remark}[section]

\begin{document}

\title[Standing waves of NLSE/NLDE on graphs with localized nonlinearity]{An overview on the standing waves of nonlinear Schr\"odinger and Dirac equations on metric graphs with localized nonlinearity}

\author[W. Borrelli]{William Borrelli}
\address{Universit\'e Paris-Dauphine, PSL Research University, CNRS, UMR 7534, CEREMADE, F-75016 Paris, France.}
\email{borrelli@ceremade.dauphine.fr}
\author[R. Carlone]{Raffaele Carlone}
\address{Universit\`{a} ``Federico II'' di Napoli, Dipartimento di Matematica e Applicazioni ``R. Caccioppoli'', MSA, via Cinthia, I-80126, Napoli, Italy.}
\email{raffaele.carlone@unina.it}
\author[L. Tentarelli]{Lorenzo Tentarelli}
\address{Universit\`{a} ``Federico II'' di Napoli, Dipartimento di Matematica e Applicazioni ``R. Caccioppoli'', MSA, via Cinthia, I-80126, Napoli, Italy.}
\email{lorenzo.tentarelli@unina.it}

\date{\today}

\begin{abstract} 
 We present a brief overview on the existence/nonexistence of standing waves for the NonLinear Schr\"odinger and the NonLinear Dirac Equations (NLSE/NLDE) on metric graphs with localized nonlinearity. We first focus on the NLSE, both in the subcritical and the critical case, and then on the NLDE, highlighting similarities and differences with the NLSE. Finally, we show how the two equations are related in the nonrelativistic limit, proving the convergence of bound states.
\end{abstract}

\maketitle

\vspace{-.5cm}
\smallskip

    

\section{Introduction}
\label{sec-intro}

The aim of this paper is to present the state of the art on the study of the standing waves of the NonLinear Schr\"odinger and the NonLinear Dirac Equations (NLSE/NLDE) on metric graphs with \emph{localized nonlinearities}.

In the following, by \emph{metric graph} we mean the locally compact metric space which one obtains endowing a \emph{multigraph} $\G=(\V,\E)$ with a parametrization that associates each bounded edge $e\in\E$ with a closed and bounded interval $I_e=[0,\ell_e]$ of the real line, and each unbounded edge $e\in\mathrm{E}$ with a (copy of the) half-line $I_e=\R^+$ (an extensive description can be found in \cite{AST-CVPDE,BK} and references therein). As a consequence, functions on metric graphs $u=(u_e)_{e\in\mathrm{E}}:\G\to\R,\,\C$ have to be seen as bunches of functions $u_e:I_e\to\R,\,\C$ such that $u_{|_e}=u_e$. Consistently, Lebesgue and Sobolev spaces are defined as
\[
 L^p(\G):=\bigoplus_{e\in\mathrm{E}}L^p(I_e),\quad p\in[1,\infty],\qquad\text{and}\qquad H^m(\G):=\bigoplus_{e\in\mathrm{E}}H^m(I_e),\quad m\in\N,
\]
and are equipped with the natural norms. Moreover, throughout $L^p$-norms are denoted by $\|u\|_{p,\G}$ and the $H^1$-norm is denoted by $\|u\|$, for the sake of simplicity.

\begin{figure}[t]
\centering
\includegraphics[width=.55\columnwidth]{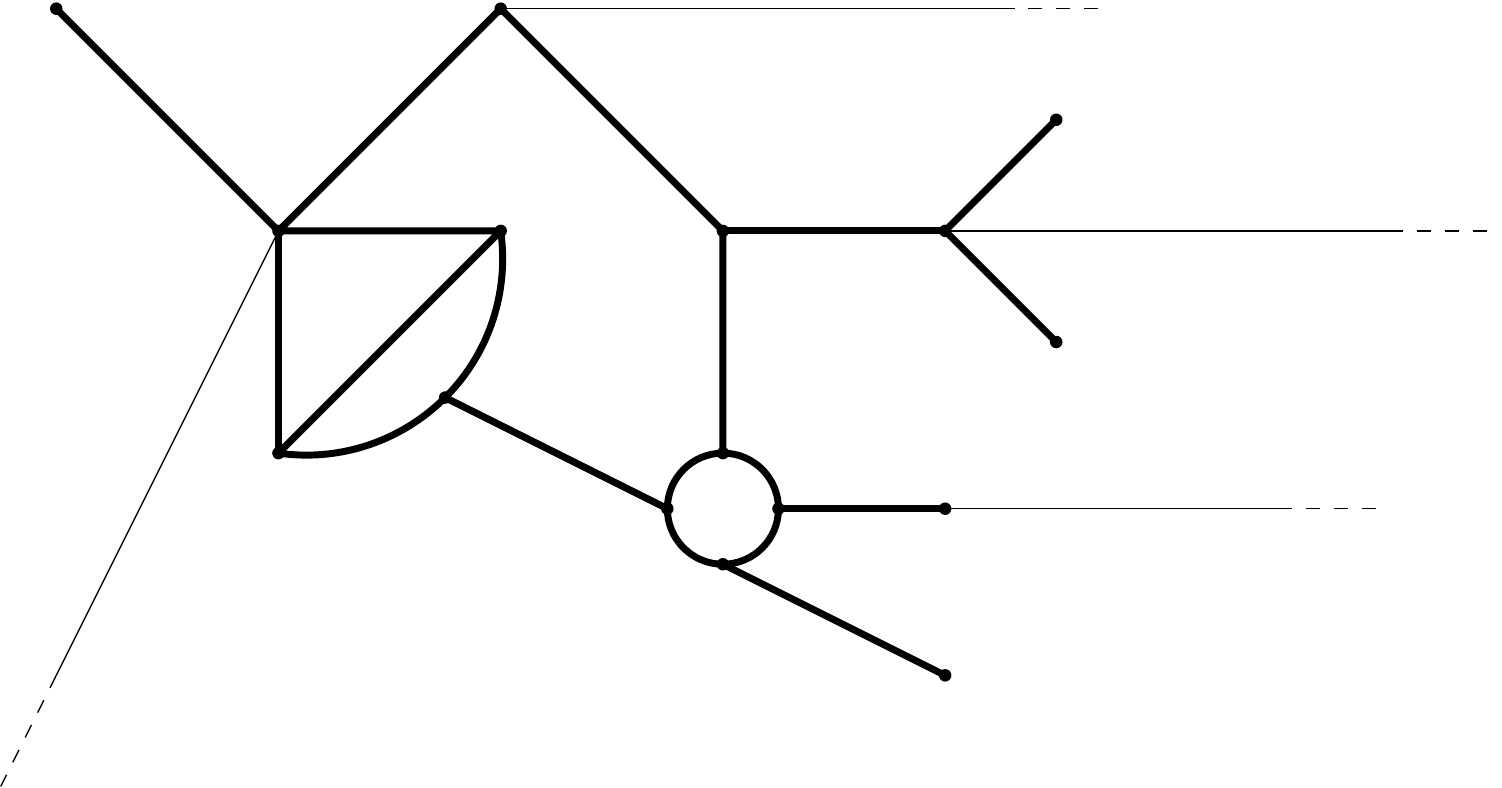}
\caption{A graph $\G$ satisfying {\bf (H1)-(H3)} and its compact core (bold edges).}
\label{fig-gen}
\end{figure}

It is also worth recalling that, in the present paper we limit ourselves to focus on the case of metric graphs $\G$ satisfying the following hypothesis:
\begin{itemize}
 \item[\textbf{(H1)}] $\G$ is \emph{connected} and \emph{noncompact};
 \item[\textbf{(H2)}] $\G$ has a \emph{finite number of edges};
 \item[\textbf{(H3)}] $\G$ has a non-empty \emph{compact core} $\K$ (which is the subgraph of $\G$ consisting of all its bounded edges)
\end{itemize}
(see, e.g., Figure \ref{fig-gen}).

\medskip
The growing interest in the study of evolution equations on \emph{metric graphs} and networks is due to the fact that they are regarded as effective models for the dynamics of physical systems constrained in branched spatial structures (see \cite{AST-PARMA} and references therein).

In particular, motivated by the study of qualitative properties of Bose-Einstein Condensates in ramified traps (see, e.g., \cite{GW-PRE}), in the last years a considerable attention has been devoted to the \emph{focusing} NLSE, i.e.
\begin{equation}
 \label{eq-NLStime}
 \imath\partial_t w=-\lap w-|w|^{p-2}\,w,\qquad p\geq2,
\end{equation}
where $-\lap$ is a suitable self-adjoint realization of the operator
\[
 -\Delta_{|\oplus_{e\in\E}C_0^\infty(\mathring{I_e})},
\]
and, precisely, on the existence of \emph{standing waves} of \eqref{eq-NLStime}. Those are solutions of the form 
$$
w(t,x)=e^{-i\lambda t}\,u(x)\,,
$$ 
with $\lambda\in\R$ and $u\in L^2(\G)$ solving the stationary version of  \eqref{eq-NLStime}, namely,
\begin{equation}
   \label{eq-NLS}
   -\lap u-|u|^{p-2}\,u=\lambda u.
\end{equation}

\begin{figure}[t]
\centering
{\includegraphics[width=.35\columnwidth]{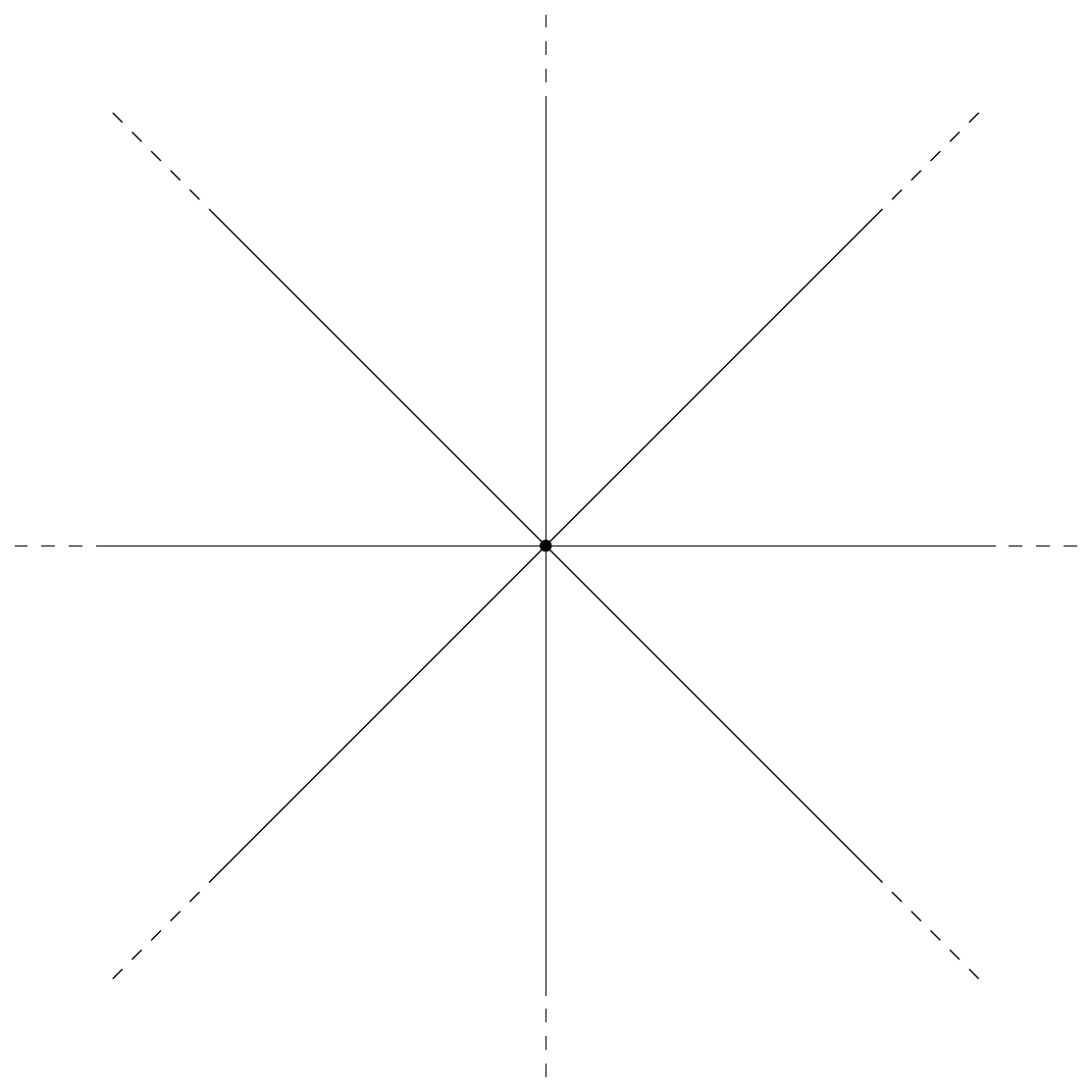}}
\caption{Infinite $N$-star graph ($N=8$).}
\label{fig-nstar}
\end{figure}

The first results in this direction (e.g., \cite{ACFN-RMP,ACFN-JPA,ACFN-JDE,ACFN-ANHIPC}) considered the so-called \emph{infinite $N$-star} graph (see, e.g., Figure \ref{fig-nstar}), in the case where $-\lap$ is the Laplacian with $\delta$\emph{-type} vertex conditions, that is
\begin{gather*}
 (-\lap^{\delta,\alpha} u)_{|I_e}:=-u_e''\qquad\forall e\in\E,\qquad\forall u\in\d(-\lap^{\delta,\alpha}),\\[.2cm]
 \d(-\lap^{\delta,\alpha}):=\left\{u\in H^2(\G):u\text{ satisfies \eqref{eq-cont} and \eqref{eq-delta}}\right\},
\end{gather*}
where
\begin{gather}
 \label{eq-cont} u_{e_1}(\v)=u_{e_2}(\v),\qquad\forall e_1,e_2\succ \v,\qquad\forall\v\in\K,\\[.2cm]
 \label{eq-delta} \sum_{e\succ v}\f{du_e}{dx_e}(\v)=\alpha u(\v),\qquad\forall\v\in\K
\end{gather}
for some $\alpha\in\R$ ($e\succ \v$ meaning that the edge $e$ is incident at the vertex $\v$, while $\tf{du_e}{dx_e}(\v)$ stands for $u_e'(0)$ or $-u_e'(\ell_e)$ according to whether $x_e$ is equal to $0$ or $\ell_e$ at $\v$).

On the other hand, in the case $\alpha=0$ which is usually called \emph{Kirchhoff} Laplacian (and which will be denoted simply by $-\lap$ in place of $-\lap^{\delta,0}$ in the sequel, for the sake of simplicity) more general graphs have been studied (precisely, any graph satisfying {\bf (H1)-(H2)}). We mention, in this regard, \cite{AST-CVPDE,AST-JFA,AST-CMP} for a discussion of the existence of \emph{ground states} (i.e., those standing waves that minimize the energy functional associated with \eqref{eq-NLS}) and \cite{AST-CVPDE2,CFN-PRE,KP-JDE,NPS-Non,NRS-JDE} concerning, more generally, \emph{excited states}. We also mention \cite{CFN-Non} where the same problems are studied in the presence of an external potential.

\medskip
A modification of this model, introduced in \cite{GSD-PRA,N-RSTA}, consists of assuming that the nonlinearity affects only the compact core of the graph (which then must be supposed non empty as in {\bf (H3)}), so that \eqref{eq-NLS} reads
\begin{equation}
   \label{eq-NLSconc}
   -\lap u-\chi_{\K}|u|^{p-2}\,u=\lambda u.
\end{equation}
where $\chi_{\K}$ is the characteristic function of $\K$. The existence of stationary solutions to \eqref{eq-NLSconc} has been discussed in \cite{ST-JDE,ST-NA,T-JMAA} in the $L^2$-subcritical case $p\in(2,6)$ and, more recently, in \cite{DT-OTAA,DT-p} in the critical case $p=6$. A detailed description of these results will be presented in Section \ref{sec-nls}.

\begin{figure}[t]
\centering
\subfloat[][one-dimensional periodic graph]
{\includegraphics[width=.65\columnwidth]{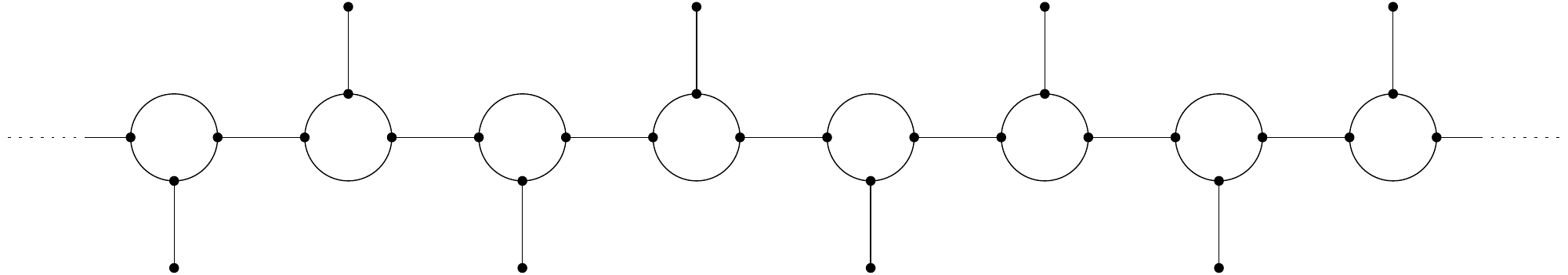}}
\quad
\subfloat[][two-dimensional grid]
{\includegraphics[width=.65\columnwidth]{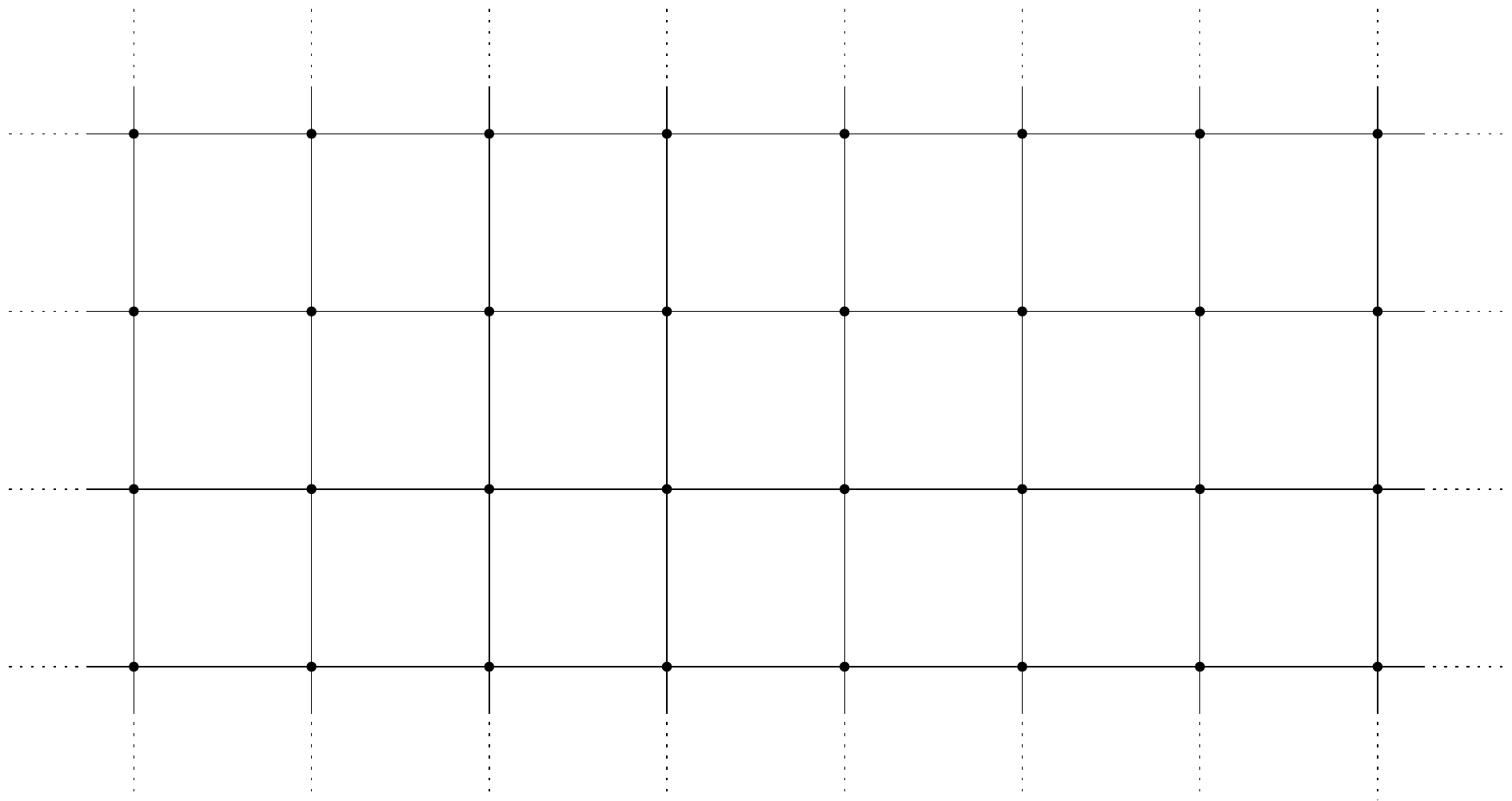}}
\caption{Two examples of periodic graphs.}
\label{fig-period}
\end{figure}

For the sake of completeness, we also remark that the NLSE on compact graphs (which, in particular, do not fulfill {\bf (H1)}) has been studied, e.g., in  \cite{CDS-MJM,D-JDE,D-p,MP-AMRX};
while the case of one or higher-dimensional \emph{periodic} graphs (which, in particular, do not fulfill {\bf (H2)} as, for instance, in Figure \ref{fig-period}) has been addressed, e.g., by \cite{ADST-APDE,DD-p,GPS-NODEA,PS-ANIHPA}.

Besides the NLS, other dispersive PDEs on metric graphs have been explored in the last years. We mention, for instance, the case of \cite{MNS-APDE} which deals with the dynamics for the Airy equation (motivated by the study of the KDV equation) on star graphs. 
 
\medskip
 
Recently, we started a new research project concerning the NonLinear Dirac Equation (NLDE) on metric graphs in \cite{BCT-p}. The physical motivations for such a model mainly come from solid state physics and nonlinear optics. The existence of Dirac solitons in Bose-Einstein condensates and optical lattices
and their concrete realization in discrete waveguide arrays have been recently investigated in \cite{HC-PD,TLB-AP}. In that case one may expect to recover the metric graph model in an appropriate scaling regime. We also mention that a rigorous mathematical study of the dynamics and the existence of Dirac solitons on lattices has been recently treated in \cite{AS-JMP, B-JDE,B-CVPDE,B-JMP,FW-CMP}.

\begin{figure}[t]
\centering
{\includegraphics[width=.35\columnwidth]{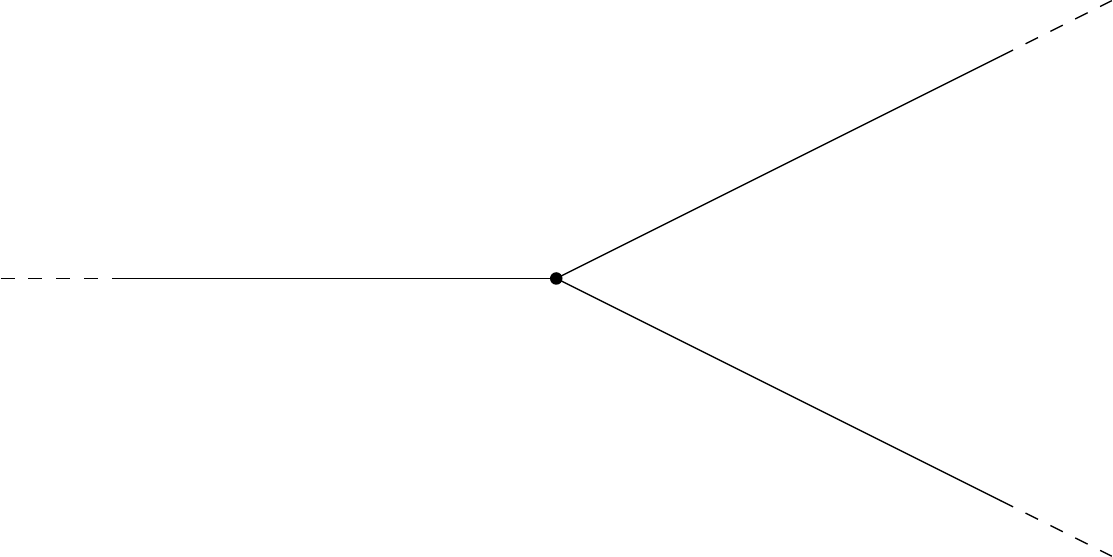}}
\caption{Infinite $3$-star graph.}
\label{fig-3star}
\end{figure}

The investigation of stationary solution of the NLDE in this context has been first proposed by \cite{SBMK-JPA} for the case of the infinite 3-star graph (see Figure \ref{fig-3star}). However, in \cite{SBMK-JPA} the authors considered the case of an \emph{extended} nonlinearity, i.e.
\[
 \imath\partial_t\Psi=\Dg\Psi-|\Psi|^{p-2}\,\Psi,\qquad p\geq2,
\]
where $\Dg$ is a suitable self-adjoint realization of the one-dimensional Dirac operator $\D$, i.e.
\[
 \D:=-\imath c\frac{d}{dx}\otimes\sigma_{1}+mc^{2}\otimes\sigma_{3},
\]
$m>0$ and $c>0$ representing the \emph{mass} of the generic particle of the system and a \emph{relativistic parameter} (respectively), and $\sigma_1$ and $\sigma_3$ representin the \emph{Pauli matrices}, i.e.
\begin{equation}
 \label{eq-pauli}
 \sigma_1:=\begin{pmatrix}
  0 & 1 \\
  1 & 0
 \end{pmatrix}
 \qquad\text{and}\qquad
 \sigma_3:=\begin{pmatrix}
  1 & 0 \\
  0 & -1
 \end{pmatrix}.
\end{equation}
Here the equation of the standing waves $\Psi(t,x)=e^{-\imath\omega t}\psi(x)$ reads
\[
 \Dg \psi-|\psi|^{p-2}\,\psi=\omega \psi.
\]
In fact, in order to deal with more complex graph topologies, in \cite{BCT-p} we restricted ourselves to the case of a Kirchhoff-type extension of the Dirac operator (for details see Section \ref{sec-D}) and, most importantly, we considered the case of a \emph{localized} nonlinearity, that is
\begin{equation}
   \label{eq-NLDconc}
   \Dg \psi-\chi_{\K}|\psi|^{p-2}\,\psi=\omega \psi.
\end{equation}
In addition, in order to test the consistency of this model we also proved (see Section \ref{sec-limit}) the convergence (of the bound states) of \eqref{eq-NLDconc} to (the bound states of) \eqref{eq-NLSconc} in the \emph{nonrelativistic limit}, namely as ``$c\to\infty$''. 

Note that in \eqref{eq-NLDconc} (as in the equation studied in \cite{SBMK-JPA}) the nonlinearity is a pure power, thus being manifestly \emph{non covariant}. This kind of nonlinearities typically arises in nonlinear optics. We also stress that from a theoretical point of view there is no conceptual contradiction as the NLDE should be intended only as an effective model.

Nevertheless, the equation obtained replacing the nonlinearity $\vert\psi|^{p-2}\,\psi$ with a \emph{covariant} one, i.e.
\[
   \Dg \psi-\vert\langle\psi,\sigma_3\psi\rangle\vert^{(p-2)/2}\psi=\omega \psi
\]
(the nonlinearity now being defined in terms of a Lorentz-scalar) is great interest. This equation might, indeed, be interpreted as the model for a \emph{truly relativistic} particle (a Dirac fermion), whose dynamics is constrained to a one-dimensional structure. In this case the model should be Lorentz-covariant (strictly speaking, in this case the equation should be rewritten in covariant form). We expect that this will require some technical adjustments in order to deal with the indefiniteness of the nonlinearity. This might be quite relevant in the study of the nonrelativistic limit, since in \cite{BCT-p} we used the fact that a positive nonlinearity allows to get a priori estimates uniform in $c$.
 
Finally, we mention that the discussion of extended nonlinearities presents some technical problems that we were not able to overcome thus far. However, we plan to investigate this question in a forthcoming paper.


\section{Nonlinear Schr\"odinger equation}
\label{sec-nls}

The study of the existence/nonexistence and multiplicity of bound states of \eqref{eq-NLSconc} carried on in \cite{DT-OTAA,DT-p,ST-JDE,ST-NA,T-JMAA} is completely based on variational methods.

This is due to the simple observation that $L^2$-solutions of \eqref{eq-NLSconc}, with $-\lap$ denoting the Kirchhoff Laplacian, arise also as \emph{constrained} critical points of the energy functional
\[
 E(u,\K,p):=\frac12\int_\G|u'|^2\dx-\frac1p\int_\K|u|^p\dx,\qquad p\geq2,
\]
on the manifold
\begin{equation}
 \label{eq-manifold}
 \Hm:=\left\{u\in\H:\|u\|_{2,\G}^2=\mu\right\}
\end{equation}
for some fixed (but arbitrary) $\mu>0$, where
\[
 \H:=\left\{u\in H^1(\G):u\text{ saisfies \eqref{eq-cont}}\right\}.
\]
Moreover, if $u$ is a constrained critical point, also known as \emph{bound state}, then an easy computation shows that the Lagrange multiplier $\lambda$ can be found as
\begin{equation}
 \label{eq-lag}
 \lambda=\lambda(u):=\frac{1}{\mu}\bigg(\int_\G|u'|^2\dx-\int_\K|u|^p\dx\bigg).
\end{equation}

\begin{rem}
 Note that $\H$ is the form domain of $-\lap$. In addition, it is worth mentioning that the choice of the manifold \eqref{eq-manifold} is also related to the fact that the time dependent counterpart of \eqref{eq-NLSconc} is mass-preserving.
\end{rem}
   

\subsection{Ground states}

The first step in the study of the bound states consists in looking for \emph{ground states}, namely the constrained minimizers of $E(\cdot,\K,p)$ on the manifold $\Hm$.

Ground states can be shown (see \cite{AST-CVPDE}) to be of the form 
$$
u(x)=e^{\imath\phi}v(x)\,
$$ 
where $\phi$ is a fixed phase factor and $v$ is real valued. Hence, in the following, in minimization problems we will always limit ourselves to consider real-valued functions.

For more general bound states we cannot prove at the moment that such a property holds. Therefore, while in existence and multiplicity results we will limit to real-valued functions, nonexistence results concern also complex-valued functions.


\subsubsection{The subcritical case: $p\in(2,6)$}

We start reporting on the existence of constrained minimizers in the so-called $L^2$-subcritical case, namely, when $p\in(2,6)$.

For $p>2$ and for every graph satisfying {\bf (H1)-(H2)}, it is possible to prove (see, e.g., \cite[Proposition 2.1]{AST-JFA} and \cite[Proposition 4.1]{T-JMAA} for real-valued functions and \cite[Proposition 2.6]{ST-NA} for complex-valued functions) that the following \emph{Gagliardo-Nirenberg inequalities} hold:
\begin{gather}
 \label{eq-GNp} \|u\|_{p,\G}^p\leq \Cp\|u'\|_{2,\G}^{\f{p}{2}-1}\|u\|_{2,\G}^{\f{p}{2}+1},\qquad\forall u\in\H,\\[.2cm]
 \label{eq-GNi} \|u\|_{\infty,\G}\leq \Ci\|u'\|_{2,\G}^{\f{1}{2}}\|u\|_{2,\G}^{\f{1}{2}},\qquad\forall u\in\H\,
\end{gather}
with $\Cp$ and $\Ci$ denoting the optimal constants of the respective inequalities. Then
\[
 E(u,\K,p)\geq\frac12\|u'\|_{2,\G}^2-\frac{\Cp}{p}\|u'\|_{2,\G}^{\f{p}{2}-1}\|u\|_{2,\G}^{\f{p}{2}+1},\qquad\forall u\in\H
\]
so that
\begin{equation}
 \label{eq-bound}
 E(u,\K,p)\geq \f{1}{2}\|u'\|_{2,\G}^{\f{p}{2}-1}\left(\|u'\|_{2,\G}^{\f{6-p}{2}}-\f{2\mu^{\f{p+2}{4}}\Cp}{p}\right),\qquad\forall u\in\Hm
\end{equation}
and, if $p\in(2,6)$, this immediately entails that the functional $E(\cdot,\K,p)$ is bounded from below on the manifold $\Hm$.

\begin{rem}
 Note that, in fact, \eqref{eq-GNp} and \eqref{eq-GNi} hold for every $u\in H^1(\G)$ under the sole assumption {\bf (H2)}, up to a redefinition of the optimal constants. However, we chose to mention the $\H$-version, which holds even if {\bf (H2)} is not fullfilled. We remark that the best constants of this version play a crucial role in the analysis of the ground states.
\end{rem}

In addition to \eqref{eq-bound}, as a general fact, one can prove that
\begin{equation}
 \label{eq-inf1}
 \inf_{u\in\Hm}E(u,\K,p)\leq0
\end{equation}
and that
\begin{equation}
 \label{eq-inf2}
 \inf_{u\in\Hm}E(u,\K,p)<0\qquad\Longrightarrow\qquad\text{a ground state does exist}.
\end{equation}
As a consequence, the following results can be established

\begin{thm}[{\cite[Theorems 3.3$\&$3.4]{T-JMAA},\cite[Corollary 3.4]{ST-NA}}]
 \label{thm-ground}
 Let $\G$ satisfy {\bf (H1)-(H3)} and let $p\in(2,6)$. Therefore:
 \begin{itemize}
  \item[1)] if $p<4$, then there exists a ground state for every $\mu>0$;\\[-.2cm]
  \item[2)] if $p\geq4$, then:\\[-.4cm]
  \begin{itemize}
   \item[(i)] whenever 
   \begin{equation}
    \label{eq-cond1}
    \mu^{\f{p-2}{6-p}}|\K|>N^{\f{4}{6-p}}c_p,
   \end{equation}
   where $N$ is the number of half-lines of $\G$ and
   \[
    c_p=\bigg[\bigg(\f{p(p-4)}{16}\bigg)^{\f{2}{p-2}}+\f{p}{8}\bigg(\f{p(p-4)}{16}\bigg)^{\f{4-p}{p-2}}\bigg]^{\f{p-2}{6-p}},
   \]
   there exists a ground state of mass $\mu$;
   \item[(ii)] whenever
   \begin{equation}
   \label{eq-cond2}
    \mu^{\f{p-2}{6-p}}|\K|<\bigg(\f{p}{2}\bigg)^{\f{2}{6-p}}\f{\Cp^{\f{4-p}{6-p}}}{\Ci^p},
   \end{equation}
   there does not exist any ground state of mass $\mu$.
  \end{itemize}
 \end{itemize}
\end{thm}

The proof of the previous theorem is clearly based on \eqref{eq-inf1} and \eqref{eq-inf2}. Precisely, from \eqref{eq-inf2} there results that the existence part can be obtained exhibiting a function with negative energy, usually called \emph{competitor}. This can be achieved for every mass $\mu$ when $p\in(2,4)$, while when $p\in[4,6)$ this is possible only under \eqref{eq-cond1}. Such competitors are constant on the compact core and exponentially decreasing on half-lines. In this way they minimize the kinetic energy on the compact core and have the exact qualitative behavior of the bound states on the half-lines (where the problem is linear). Clearly, such competitors cannot be minimizers since they do not fulfill \eqref{eq-delta} with $\alpha=0$.

On the other hand, \eqref{eq-inf1} shows that, in order to prove nonexistence, it is sufficient to prove that every function possesses a strictly positive energy level. This is actually, the idea behind the first proof of the nonexistence result for ground states. However, there is a different strategy, which is similar to that of the proof of Theorem \ref{thm-nonex}, which allows to prove item 2)(ii) in a more straight way and with the sharper threshold given by \eqref{eq-cond2}.

Finally, it is worth mentioning that the negativity of the energy levels of the minimizers entails that the associated Lagrange multiplier, given by \eqref{eq-lag}, is negative and, especially, that the results of Theorem \ref{thm-ground} are invariant under the \emph{homotetic} transformation
\[
 \mu\mapsto\sigma\mu\qquad\G\mapsto\sigma^{\f{2-p}{6-p}}\G.
\]
In other words, if one sets for instance $p=4$, a problem on a graph $\G$ with a given mass constraint $\mu$ is completely equivalent to a problem with half the mass and a ``doubled'' graph, and viceversa.


\subsubsection{The critical case: $p=6$}

The study of the so-called $L^2$-critical case $p=6$ has been successfully managed only very recently. Here the main difficulty comes from the fact that the boundedness from below of the energy functional strongly depends on the the mass constraint.

Indeed, in this case \eqref{eq-bound} reads
\begin{equation}
\label{eq-bound_bis}
E(u,\K,6)\geq \f{1}{2}\|u'\|_{2,\G}^{2}\left(1-\f{\mu^{2}\Cs}{3}\right),\qquad\forall u\in\Hm
\end{equation}
and hence the boundedness of functional clearly depends on the value of $\mu$. In particular, it turns out that the discriminating value is related to the best constant of \eqref{eq-GNp}, that is, the \emph{critical mass} $\mu_\G$ is defined by
\[
 \mu_\G:=\sqrt{\f{3}{\Cs}}\, ,
\]
and is well known when $\G=\R,\R^+$, i.e.
\[
 \mu_\R=\f{\pi\sqrt{3}}{2},\qquad\mu_\R^+=\f{\mu_\R}{2}.
\]
In addition, there holds
\[
 \mu_\R^+\leq\mu_\G\leq\mu_\R,\qquad\forall\G\text{ satisfying {\bf (H1)-(H3)}},
\]
and $\mu_\R,\mu_\R^+$ are the sole values of the mass constraint at which the problem with the extended nonlinearity on $\R,\R^+$ (respectively) admits a ground state (see, e.g., \cite{C}).

However, one can easily see that, while for the problem with extended nonlinearity (treated in \cite{AST-CMP}) $\mu_\G$ is the proper parameter to be investigated, in the localized case one should better consider a \emph{reduced} critical mass
\[
 \mu_\K:=\sqrt{\frac{3}{\Ck}},
\]
where $\Ck$ denotes the optimal constant of the following modified Gagliardo-Nirenberg inequality
\[
 \|u\|_{6,\K}^6\leq \Ck\|u'\|_{2,\G}^{2}\|u\|_{2,\G}^{4},\qquad\forall u\in\H,
\]
namely
\[
 \Ck:=\sup_{u\in \H}\f{\|u\|_{6,\K}^6}{\|u'\|_{2,\G}^{2}\|u\|_{2,\G}^{4}}.
\]
In fact, using this new parameter one obtains
\[
E(u,\K,6)\geq \f{1}{2}\|u'\|_{2,\G}^{2}\left(1-\f{\mu^{2}\Ck}{3}\right),\qquad\forall u\in\Hm,
\]
which is clearly sharper than \eqref{eq-bound_bis}, as one can easily see that $\Ck\leq \Cs$.

\begin{figure}[t]
\centering
{\includegraphics[width=.45\columnwidth]{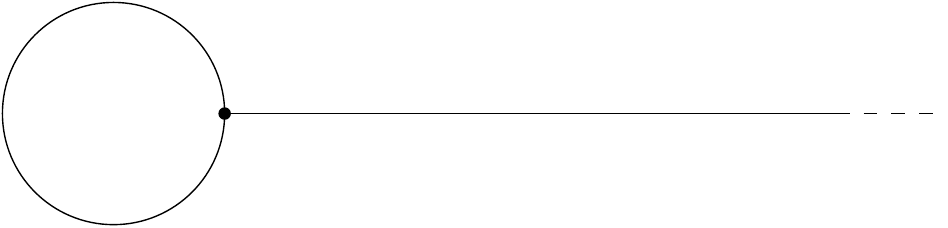}}
\caption{Tadpole graph.}
\label{fig-tadpole}
\end{figure}

\medskip
Existence of ground states has been first investigated in \cite{DT-OTAA} in the simple case of the \emph{tadpole} graph (see Figure \ref{fig-tadpole}). Very recently, more refined results have been obtained without any assumption on the topology of the graph, thus providing a (almost) complete classification of the phenomenology.

In order to state such results, let
\[
 \mathcal{E}_\G(\mu):=\inf_{u\in\Hm}E(u,\K,6),
\]
and recall the following definitions for a metric graph:
\begin{itemize}
 \item[-] a graph $\G$ is said to admit a \emph{cycle covering} if and only if every edge of $\G$ belongs to a \emph{cycle}, namely either a \emph{loop} (i.e., a closed path of consecutive bounded edges) or an unbounded path joining the endpoints of two distinct half-lines (which are then identified as a single vertex at infinity);
 \item[-] a graph $\G$ is said to possesses a \emph{terminal edge} if and only if it contains an edge ending with a vertex of degree one.
\end{itemize} 

\begin{figure}[t]
\centering
\subfloat[][a graph with one terminal edge]
{\includegraphics[width=.35\columnwidth]{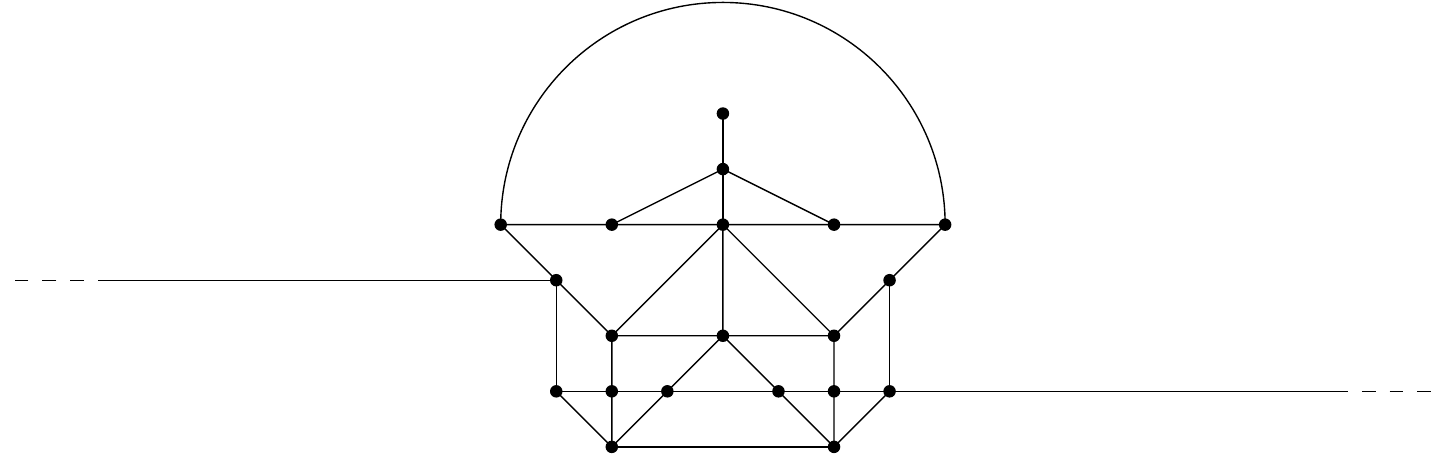}}\qquad
\subfloat[][a graph admitting a cycle-covering]
{\includegraphics[width=.35\columnwidth]{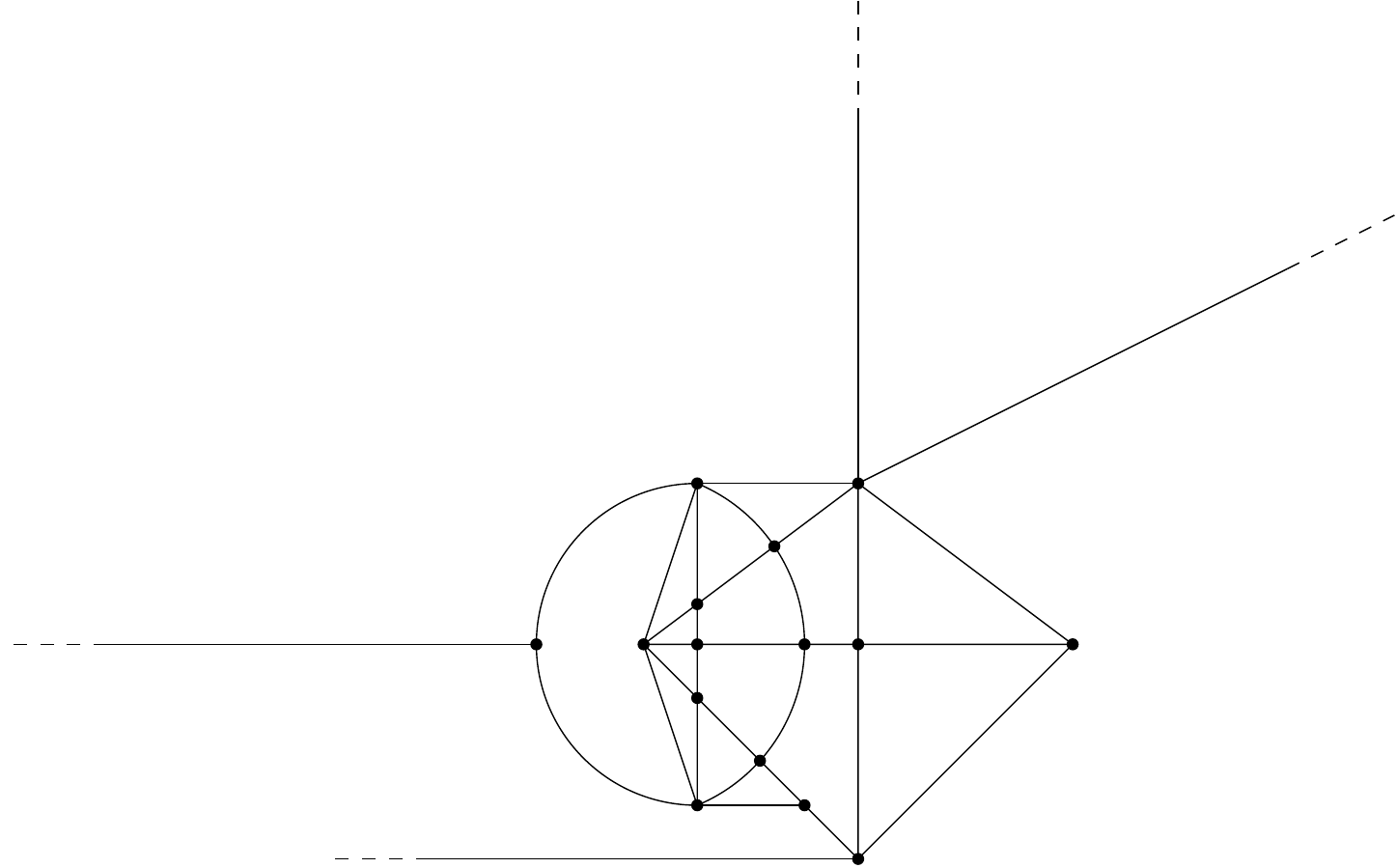}}

\medskip
\subfloat[][a graph with one half-line and without terminal edges]
{\includegraphics[width=.35\textwidth]{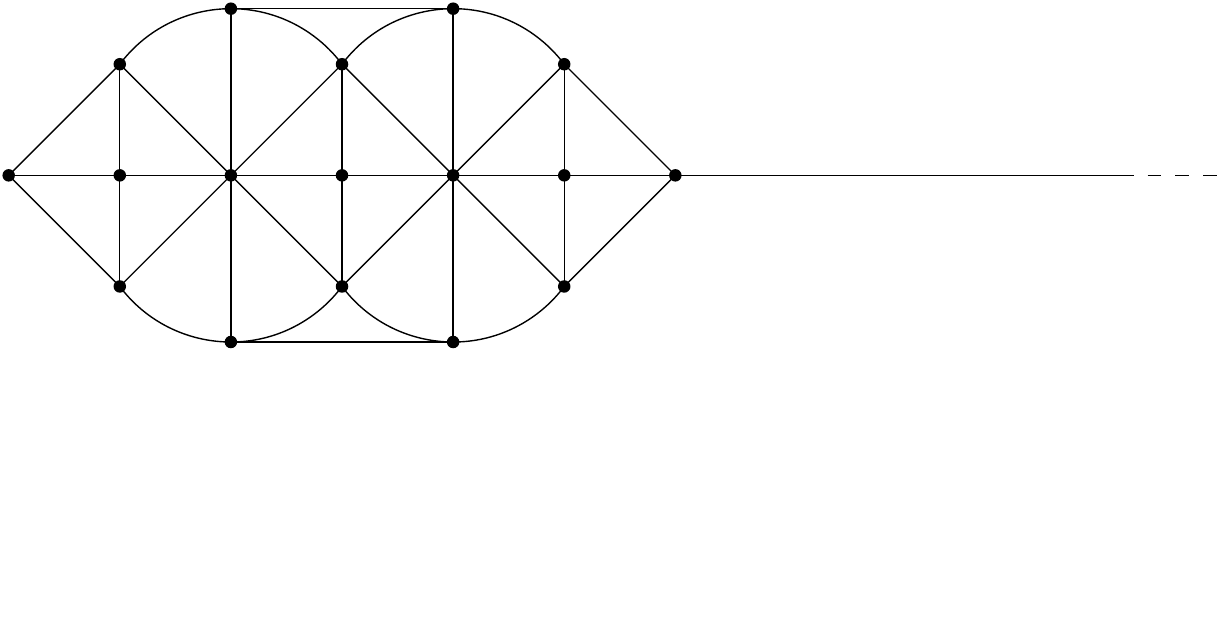}}\qquad
\subfloat[][a graph without terminal edges and cycle-coverings, and with two half-lines]
{\includegraphics[width=.35\columnwidth]{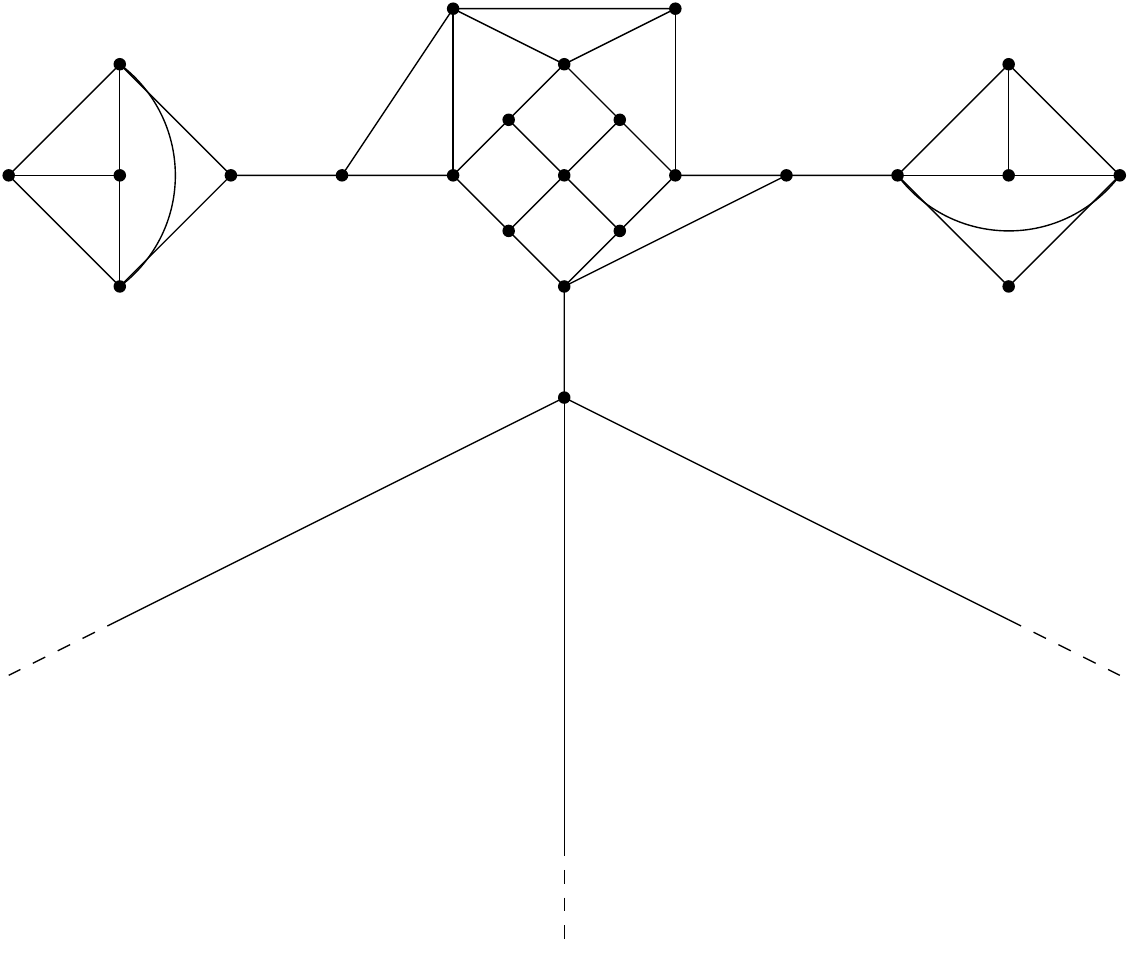}}
\caption{Examples of graphs from cases \emph{(i)-(iv)} of Theorem \ref{thm-crit}.}
\label{fig-crit}
\end{figure}

\begin{thm}[{\cite[Theorem 1.1]{DT-p}}]
 \label{thm-crit}
 Let $\G$ satisfy {\bf (H1)-(H3)}. Then,
 \[
   \mathcal{E}_\G(\mu,\K)\begin{cases}
   =0 & \text{if }\mu\leq\mu_\K\\
   <0 & \text{if }\mu\in(\mu_\K,\mu_\R]\\
   =-\infty & \text{if }\mu>\mu_\R.
  \end{cases}
 \]
 In addition,
 \begin{itemize}
 \item[(i)] if $\G$ has at least a terminal edge (as, for instance, in Figure \ref{fig-crit}{\footnotesize (A)}), then
 \[
  \mu_\K=\mu_{\R^+},\qquad \mathcal{E}_\G(\mu,\K)=-\infty\quad\text{for all }\mu>\mu_\K,
 \]
 and there is no ground state of mass $\mu$ for any $\mu>0$;
 \item[(ii)] if $\G$ admits a cycle-covering (as, for instance, in Figure \ref{fig-crit}{\footnotesize (B)}), then
 \[
  \mu_\K=\mu_\R
 \]
 and there is no ground state of mass $\mu$ for any $\mu>0$;
 \item[(iii)] if $\G$ has only one half-line and no terminal edges (as, for instance, in Figure \ref{fig-crit}{\footnotesize (C)}), then 
 \begin{equation}
  \label{eq-crit_half}
  \mu_{\R^+}<\mu_\K<\sqrt{3}
 \end{equation}
 and there is a ground states of mass $\mu$ if and only if $\mu\in[\mu_\K,\mu_\R]$.
 \item[(iv)] if $\G$ has no terminal edges, does not admit a cycle-covering, but presents at least two half-lines (as, for instance, in Figure \ref{fig-crit}{\footnotesize (D)}), then
 \begin{equation}
  \label{eq-crit_twohalf}
  \mu_{\R^+}<\mu_\K\leq\mu_\R
  \end{equation}
  and there is a ground states of mass $\mu$ if and only if $\mu\in[\mu_\K,\mu_\R]$, provided that $\mu_\K\neq\mu_\R$.
 \end{itemize}
\end{thm}

Theorem \ref{thm-crit} presents several differences with respect to its analogous for the extended case established in \cite{AST-CMP} (in addition to the fact that $\mu_\G$ is replaced with $\mu_\K$). For details on such differences we refer the reader to \cite{DT-p}. Here, we aim at highlighting a more crucial feature of the previous result. Indeed, while in the subcritical case the topology of the graph plays no role in establishing existence/nonexistence of ground states, here the topological classification of the graphs is the central point, as it occurs in the context of extended nonlinearities. Even though there is no rigorous explanation of such a phenomenon, this seems to be related to the fact that the critical power of the nonlinearity ``breaks'' the characteristic lenght of the graph, i.e. $\mu^{\f{p-2}{6-p}}|\K|$, thus making metric properties less relevant in the discussion.
 
Nevertheless, they preserve a role also in the critical case, at least in cases (iii) and (iv), as shown by the following

\begin{thm}[{\cite[Theorems 1.2$\&$1.3]{DT-p}}]
 \label{thm-asymp}
 Estimates \eqref{eq-crit_half} and \eqref{eq-crit_twohalf} are sharp in general; i.e., for every $\varepsilon>0$ there exist two non-compact metric graphs $\G_\varepsilon^1,\G_\varepsilon^2$ (with compact cores $\K^1_\ep,\K^2_\ep$), with one half-line and without terminal edges, such that 
 \[
  \mu_{\K_\varepsilon^1}\leq \mu_{\R^+}+\varepsilon\qquad\text{and}\qquad
  \mu_{\K_\varepsilon^2}\geq\sqrt{3}-\varepsilon,
 \]
 and two non-compact metric graphs $\G_\varepsilon^3,\G_\varepsilon^4$ (with compact cores $\K^3_\ep,\K^4_\ep$), without terminal edges and cycle coverings and with at least two half-lines, such that
 \[
  \mu_{\K_\varepsilon^3}\leq\mu_{\R^+}+\varepsilon\qquad\text{and}\qquad\mu_{\K_\varepsilon^4}\geq\mu_\R-\varepsilon.
 \]
\end{thm}

\begin{figure}[t]
\centering
\subfloat[][a graph whose compact core does not admit a cycle covering]
{\includegraphics[width=.4\columnwidth]{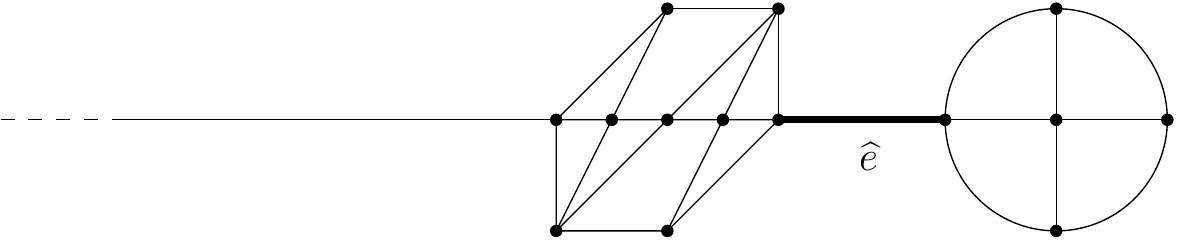}}\qquad\qquad
\subfloat[][a tadpole graph with extra connections between two vertices]
{\includegraphics[width=.4\columnwidth]{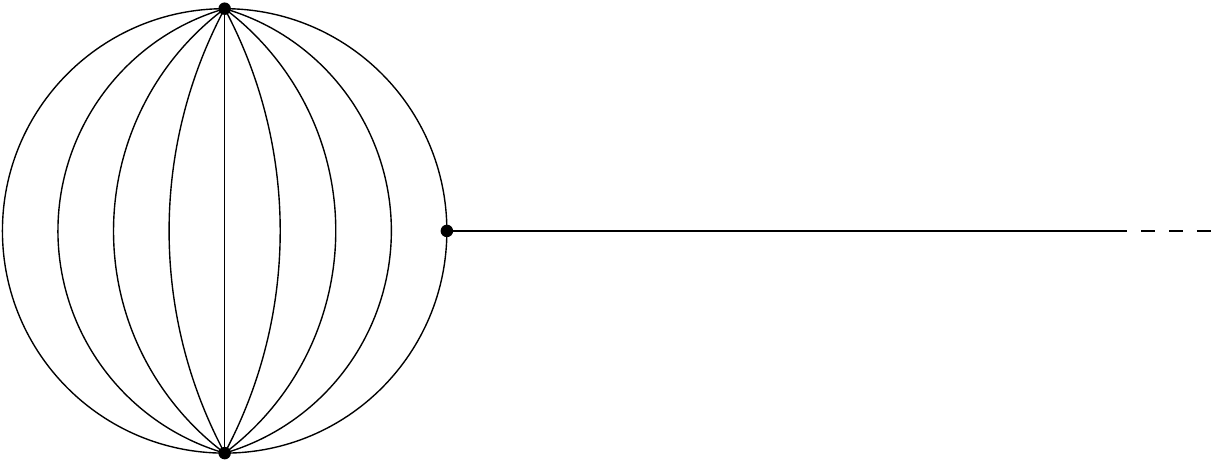}}

\medskip
\medskip
\subfloat[][a signpost graph]
{\includegraphics[width=.4\textwidth]{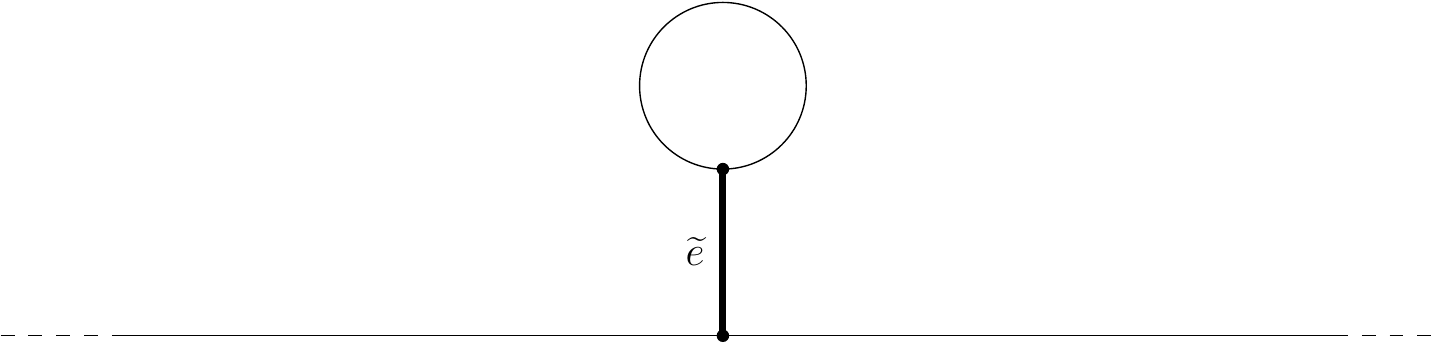}}\qquad
\caption{Examples of graphs from Theorem \ref{thm-asymp}.}
\label{fig-examples}
\end{figure}

In other words Theorem \ref{thm-asymp} shows that estimates \eqref{eq-crit_half} and \eqref{eq-crit_twohalf} are sharp by exhibiting four suitable sequences of graphs. They can be constructed as follows:
\begin{itemize}
 \item[1)] the sequence $\G_\ep^1$ can be constructed by considering a graph whose compact core does not admit a cycle covering (see, e.g., Figure \ref{fig-examples}{\footnotesize (A)}) and letting the lenght of one of its \emph{cut-edges}, the edges whose removal disconnects the graph (e.g., $\widehat{e}$ in Figure \ref{fig-examples}{\footnotesize (A)}), go to infinity;
 \item[2)] the sequence $\G_\ep^2$ can be constructed by considering a graph as in Figure \ref{fig-examples}{\footnotesize (B)}) and letting the lenght of the compact core go to infinity, keeping the total diameter of the compact core bounded (namely, thickening the compact core);
 \item[3)] the sequences $\G_\ep^3,\G_\ep^4$ can be constructed by considering a \emph{signpost} graph (see, e.g., Figure \ref{fig-examples}{\footnotesize (C)}) and letting the lenght of its cut-edge $\widetilde{e}$ go to infinity and to zero, respectively.
\end{itemize}

The proofs of Theorems \ref{thm-crit} and \ref{thm-asymp} are quite technical and we refer the reader to \cite{DT-p} for a complete presentation. They are based on a careful analysis of the behavior of the constant $\Ck$ and on the ``graph surgery'' and rearrangement techniques developed by \cite{AST-CVPDE,AST-JFA,AST-CMP}. In particular, given a nonnegative function $u\in\H$, one uses the \emph{decreasing} rearrangement
\[
 u^*(x):=\inf\{t\geq0:\rho(t)\leq x\}\quad x\in[0,|\G|)
\]
and the \emph{symmetric} rearrangement
\[
 \widehat{u}(x):=\inf\{t\geq0:\rho(t)\leq 2|x|\}\qquad x\in(-|\G|/2,|\G|/2),
\]
where
\[
 \rho(t):=\sum_{e\in\mathrm{E}}|\{x_e\in I_e:u_e(x_e)>t\}|\quad t\geq0,
\]
and precisely the facts that
\[
 \|u\|_{p,\G}=\|u^*\|_{p,\R^+}=\|\widehat{u}\|_{p,\R},\qquad\forall\,p\geq1,
\]
that
\[
 \|(u^*)'\|_{2,\R^+}\leq\|u'\|_{2,\G}
\]
and that, for every function $u$ which possesses at least two preimages at each level,
\[
 \|(\widehat{u})'\|_{2,\R}\leq\|u'\|_{2,\G},
\]
in order to construct suitable minimizing sequences or suitable competitors to exploit level arguments analogous to \eqref{eq-inf2} to get compactness in the limit.


\subsection{Bound states}

Besides the investigation of ground states, one can also study standing waves which are not necessarily minimizers of the constrained energy, namely \emph{bound/excited states}. They may arise as constrained critical points of the energy functional on $\Hm$.

To the best of our knowledge all the known results on existence (and multiplicity), as well as on nonexistence, of general bound states strongly exploit the assumption on the subcriticality of the power nonlinearity (the critical case is open thus far).


\subsubsection{Existence results}

Existence and multiplicity results in this context have been proved by extending some well known techniques from Critical Point Theory (see, e.g., \cite{AM,BL-ARMA,R-CBMS}).

However, the context of metric graphs presents two additional technical issues. First, it is not possible, in general, to gain compactness restricting to symmetric (in some suitable sense) functions if one does not want to restrict the discussion to symmetric graphs. Furthermore, the fact that graphs with non empty compact core are not homothetically invariant entails that multiple solutions of prescribed mass cannot be found by scaling arguments, in general.

The strategy used in \cite{ST-JDE}, in order to solve these problems is, then, the following:
\begin{itemize}
 \item[-] detect the energy levels at which the \emph{Palais-Smale condition} is satisfied, namely detect the values $c\in\R$ such that any sequence $(u_n)\in\Hm$ satisfying\\[-.3cm]
 \begin{itemize}
  \item[(i)] $E(u_n,\K,p) \to c$,\\[-.3cm]
  \item[(ii)] $\|dE_{|\Hm}(u_n,\K,p)\|_{\mathbb{T}'_{u_n}\Hm}\to 0$\\[-.3cm]
 \end{itemize}
 (with $\mathbb{T}'_{u_n}\Hm$ denoting the topological dual of the tangent to the manifold $\Hm$ at $u_n$), admits a subsequence converging in $\Hm$.
 \item[-] constructing suitable \emph{min-max} levels.
\end{itemize}
In this case one can prove that the constrained functional possesses the Palais-Smale property only at negative levels and whence the min-max levels have to be negative. In addition, as the functional is even, it is possible to construct such min-max levels using the \emph{Krasnosel'skii genus} (see, e.g., \cite{K}); that is, for every $A\subset \H\backslash\{0\}$ closed and symmetric, the natural number defined by
\[
\gamma(A):=\min\{n\in\N:\exists\phi:A\rightarrow\R^n\backslash\{0\},\, \phi\mbox{ continuous and odd}\}.
\]
In fact, one can prove that the suitable levels are
\[
c_j := \inf_{A\in \Gamma_j} \max_{u\in A} E(u,\K,p),
\]
where
\[
\Gamma_j := \{A \subset \Hm : A \text{ is compact and symmetric, and } \gamma(A) \ge j\}.
\]

\medskip
In this way it is possible to prove the following

\begin{thm}[{\cite[Theorem 1.2]{ST-JDE}}]
Let $\G$ satisfy {\bf (H1)-(H3)} and let $p\in(2,6)$. For every $k\in\N$, there exists ${\mu_k}>0$ such that for every $\mu \ge {\mu_k}$ there exist at least $k$ distinct pairs $(\pm u_j)$ of bound states of mass $\mu$. Moreover, for every $j=1, \dots, k$, 
\[
E(u_j,\K,p) \le jE(\varphi_{\mu/j},\R,p) + \sigma_k(\mu) < 0,
\]
where $\varphi_{\mu/j}$ denote the unique positive minimizer of $E(\cdot,\R,p)$ constrained on $\mathcal{H}_{\mu/j}(\G)$ and $\sigma_k(\mu) \to 0$ (exponentially fast) as $\mu \to \infty$. Finally, for each $j$, the Lagrange multiplier $\lambda_j$ relative to $u_j$ is negative.
\end{thm}

The strategy hinted before cannot be easily adapted to the extended problem since in this case the Palais-Smale condition fails also at infinitely many negative energy levels and then it is open how to fix a min-max level in the gap between two of these levels.

On the other hand, it is also unclear how to extend this strategy to the critical case. Indeed the construction of negative energy min-max levels is based on the possibility of putting (on increasing the mass) several suitable truncations of scaled copies of $\varphi_{\mu/j}$ on the compact core of the graph keeping the total energy negative. However, no direct analogous of such a techniques is available for the critical case as the minimizers $\varphi_\mu$ do exist only at the critical mass and present a zero energy level.


\subsubsection{Nonexistence results}

Nonexistence results concern clearly only the regime $p\in[4,6)$ since for $p\in(2,4)$ the existence of bound states is guaranteed by Theorem \ref{thm-ground}.

The interesting part of nonexistence results for bound states in the subcritical regime is that they strongly relies both on metric and on topological features of the graph as one can see by the following

\begin{figure}[t]
\centering
{\includegraphics[width=.4\columnwidth]{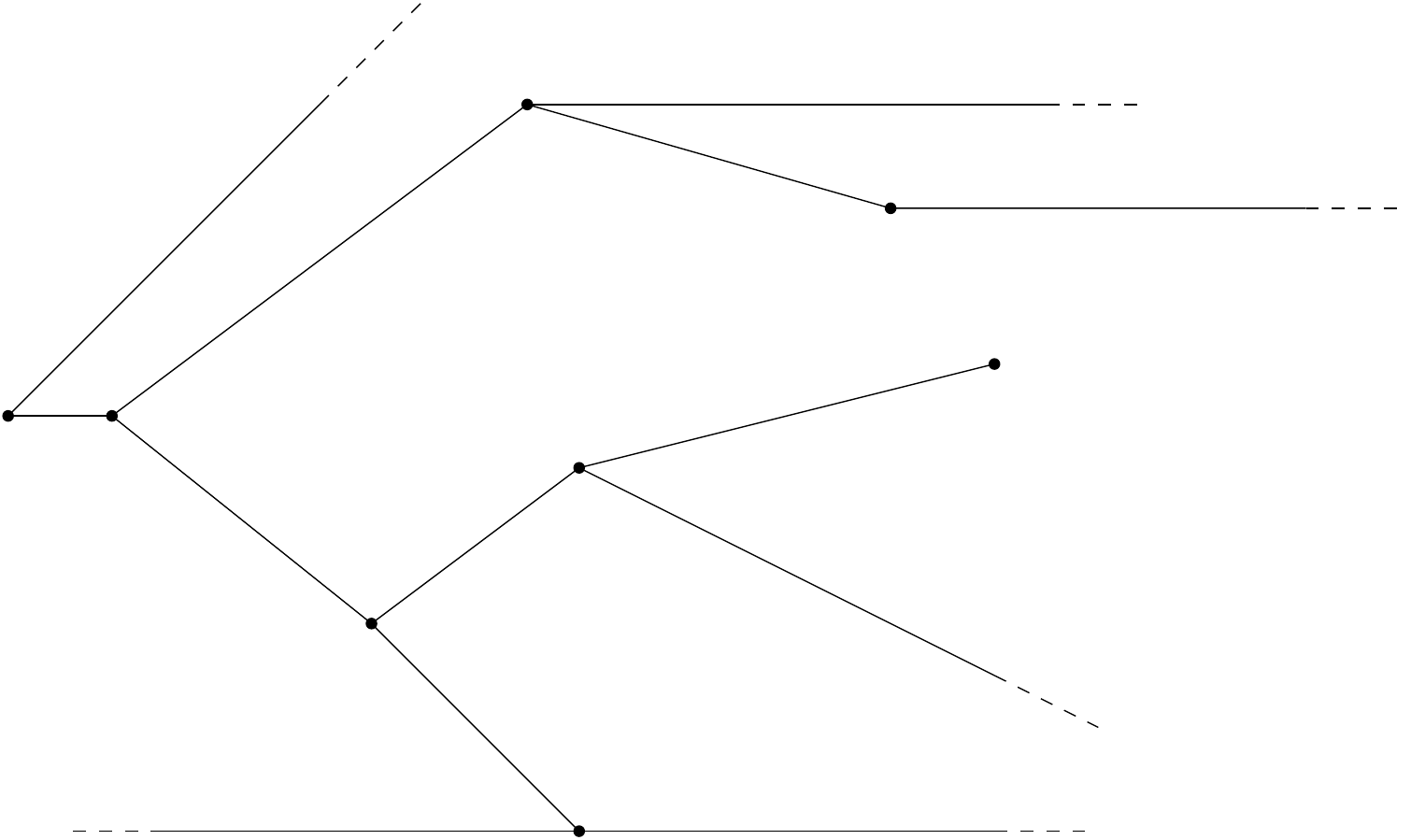}}
\caption{A tree with one pendant.}
\label{fig-tree}
\end{figure}

\begin{thm}[{\cite[Theorems 3.2$\&$3.5]{ST-NA}}]
 \label{thm-nonex}
 Let $\G$ satisfy {\bf (H1)-(H3)} and let $p\in[4,6)$. Therefore,
 \begin{itemize}
  \item[(i)] if the graph $\G$ satisfies
  \[
   \mu^{\frac{p-2}{6-p}}|\K|<\f{\Cp^{\f{4-p}{6-p}}}{\Ci^p},
  \]
  then, there are no bound states of mass $\mu$ with $\lambda\leq0$;\\[-.3cm]
  \item[(ii)] if $\G$ is a \emph{tree} (i.e., no loops) with at most one pendant (see, e.g., Figure \ref{fig-tree}), then there is no bound state of mass $\mu$ with $\lambda\geq0$, every $\mu>0$.
 \end{itemize}
\end{thm}

The most remarkable fact is that the dependence on metric or topological features in the nonexistence result is connected to the sign of the Lagrange multiplier. In particular, the first condition prevents the existence of bound states supported on the whole $\G$, as such functions cannot possess a negative Lagrange multiplier since they are in $L^2(\G)$. On the other hand, one can check that the second part of Theorem \ref{thm-nonex}, holds as well for $p\geq6$ and even in the extended nonlinearity case.

\begin{figure}[t]
\centering
\subfloat[][a graph with a loop]
{\includegraphics[width=.4\columnwidth]{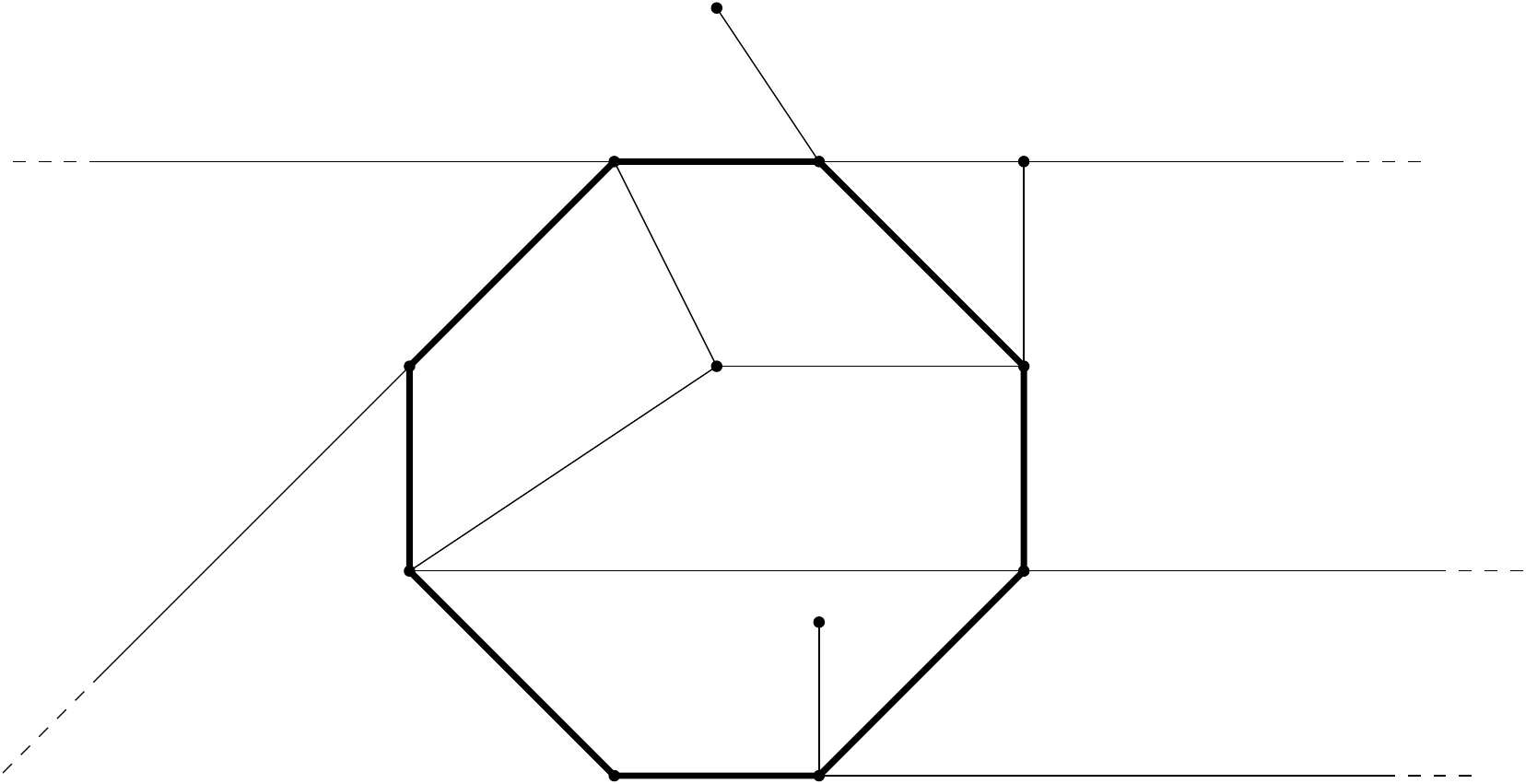}}\qquad
\subfloat[][a tree with two pendants]
{\includegraphics[width=.4\columnwidth]{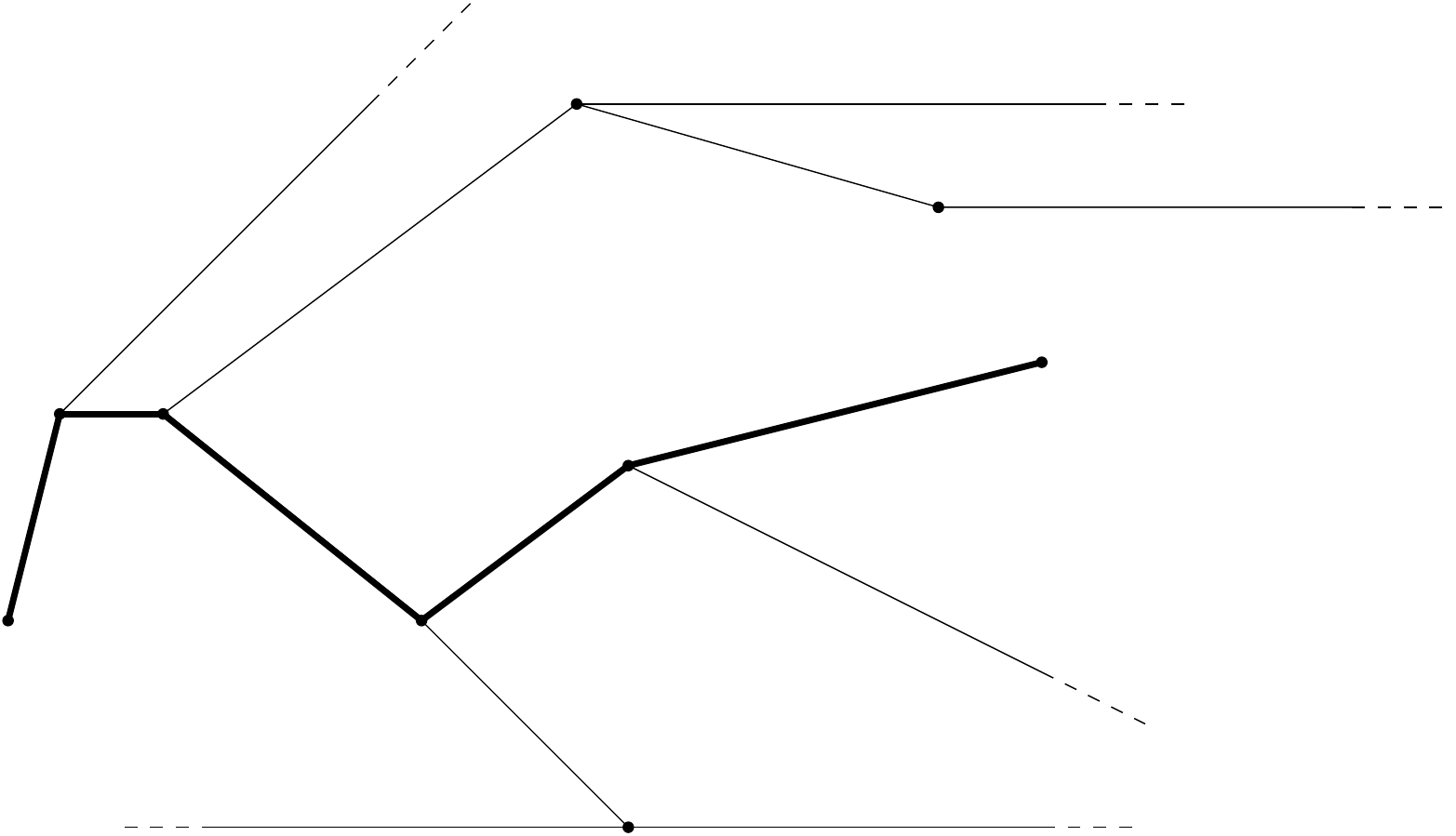}}
\caption{Examples of graphs for which item (ii) of Theorem \ref{thm-nonex} does not hold.}
\label{fig-nonoex}
\end{figure}

Finally, it is worth observing that the condition of the second part of Theorem \ref{thm-nonex} is sharp. Indeed, one can easily construct counterexamples whenever the assumption of being a tree with at most one pendant is dropped. Precisely, if $\G$ possesses a loop, then one can define a bound state with a nonnegative Lagrange multiplier supported only on the loop (see, e.g., Figure \ref{fig-nonoex}{\footnotesize (A)}); if $\G$ is a tree with two pendants, instead, then one can define a bound state with nonnegative Lagrange multiplier supported on the path which joins the two pendants (see, e.g., Figure \ref{fig-nonoex}{\footnotesize (B)}).

\section{Nonlinear Dirac equation}

As already mentioned in Section \ref{sec-intro}, the study of bound state-solutions of NLDE \eqref{eq-NLDconc} in the context of metric graphs has been introduced in \cite{SBMK-JPA}. To our knowledge, the first rigorous work on this subject for general graphs has been done in \cite{BCT-p}, for the case of localized nonlinearities.

Such discussion requires a suitable definition of the Dirac operator on graphs and the adaptation of some technques from Critical Point Theory.

Before going more into details, we recall that Dirac-type equations are of spinorial nature and, hence, here we de deal with vector-valued function on graphs (\emph{$2$-spinors}). More precisely, we consider wave functions $\psi:\G\to\C^2$, that can be seen as a family of $2$-spinors on intervals, i.e. $\psi=(\psi_e)$ with
\[
 \psi_e=\begin{pmatrix}\phi_e\\[.2cm] \chi_e\end{pmatrix}:I_e\longrightarrow\C^{2},\qquad \forall e\in\E,
\]
or equivalently as a $2$-components vector of functions on graphs, i.e. $\psi=(\phi,\chi)^T$, with $\phi=(\phi_e)$ and $\chi=(\chi_e)$. Accordingly,
\[
L^{p}(\G,\C^{2}):= \bigoplus_{e\in\E} L^{p}(I_e,\C^2)=\{\psi:\G\to\C^2:\phi,\,\chi\in L^p(\G)\},
\]
and
\[
 H^{m}(\G,\C^{2}):= \bigoplus_{e\in\E} H^{m}(I_e,\C^2)=\{\psi:\G\to\C^2:\phi,\,\chi\in H^m(\G)\},
\]
endowed with the natural norms, which we denote, with a little abuse of notation, by $\|\psi\|_{p,\G}$ and $\|\psi\|$ (in the case $m=1$), respectively.


\subsection{Remarks on the Dirac operator on graphs}
\label{sec-D}

As we mentioned before, the first step in the study of the NLDE on graphs is to define a self-adjoint realization of the Dirac operator on graphs.

A complete discussion of such a topic can be found, e.g., in \cite{BT-JMP,P}. Here, however, we limit ourselves to the extension that we call of \emph{Kirchhoff type}, since it represents the analogous to the Schr\"odinger operator with Kirchhoff conditions and, hence, corresponds to the \emph{free case}, that is, the case where there is no interaction at vertices (for detail on such extension, refer to \cite{BCT-p}).

For every fixed $m,\,c>0$, we define the Dirac operator of Kirchhoff-type as the operator 
$$
\Dg:L^2(\G,\C^2)\to L^2(\G,\C^2)
$$
 with action
\[
{\Dg}_{|I_e}\psi=\D_e\psi_e:=-\imath c\,\sigma_{1}\psi_e'+mc^{2}\,\sigma_{3}\psi_e,\qquad\forall e\in\E,\qquad\forall \psi\in\d(\Dg)
\]
$\sigma_1,\sigma_3$ being the Pauli matrices defined in \eqref{eq-pauli}, and whose domain is
\[
 \d(\Dg):=\left\{\psi\in H^1(\G,\C^2):\psi\text{ satisfies \eqref{eq-cont2} and \eqref{eq-kirch2}}\right\},
\]
where
\begin{gather}
 \label{eq-cont2}
 \phi_{e_1}(\v)=\phi_{e_2}(\v),\qquad\forall e_1,e_2\succ \v,\qquad\forall\v\in\K,\\[.2cm]
 \label{eq-kirch2}
 \sum_{e\succ \v}\chi_{e}(\v)_{\pm}=0,\qquad\forall \v\in\K,
\end{gather}
$\chi_{e}(\v)_{\pm}$ standing for $\chi_{e}(0)$ or $-\chi_{e}(\ell_e)$ according to whether $x_e$ is equal to $0$ or $\ell_e$ at $\v$.

Such operator can be proved to be self-adjoint and, in addition, possesses spectral properties analogous to the standard Dirac operator on $\R$, that is
\begin{equation}
 \label{eq-sp_D}
 \sigma(\Dg)=(-\infty,-mc^2]\cup[mc^2,+\infty)\, ,
\end{equation}
even though, according to the structure of the graph, there might be eigenvalues embedded at any point of the spectrum.

The actual reason for which we call such operator of Kirchhoff type is better explained by Section \ref{sec-limit}. However, an informal justification is provided by the following formal computation. First, one can see that $\Dg^2$ acts as $(-\lap)\otimes\mathbb{I}_{\C^2}$. In addition, if one considers spinors of the type $\psi=(\phi,0)^T$ and assumes that they belong to the domain of $\Dg^2$, namely that $\psi\in\d(\Dg)$ and that $\Dg\psi\in\d(\Dg)$, one clearly sees that $\phi\in\d(-\lap)$. 

\medskip
Finally, compared to the Schr\"odinger case, a big difference is given by the definition of the \emph{quadratic form} associated with $\Dg$. In particular, here it is not explicitly known since Fourier transform is not available (in a simple way), as we are not in an euclidean space, and since classical duality arguments fail as it is not true in general that $H^{-1/2}(\G,\C^2)$ is the topological dual of $H^{1/2}(\G,\C^2)$.

It can be defined, clearly, using the \emph{Spectral Theorem}, so that
\[
 \d(\Q):=\bigg\{\psi\in L^2(\G,\C^2):\int_{\sigma(\Dg)}|\nu|\,d\mu_{\Dg}^\psi(\nu)\bigg\},\qquad \Q(\psi):=\int_{\sigma(\Dg)}\nu\,d\mu_{\Dg}^\psi(\nu),
\]
where $\mu_{\Dg}^\psi$ denotes the spectral measure associated with $\Dg$ and $\psi$. However, such a definition is not useful for computations. A more precise description of the quadratic form and its domain can be obtained arguing as follows (for more details, see \cite{BCT-p}).

First, using Real Interpolation Theory, one can prove that
\begin{equation}
 \label{eq-domex}
 \d(\Q)=\left[L^{2}(\G,\C^{2}),\d(\Dg)\right]_{\frac{1}{2}},
\end{equation}
whence
\[
  \d(\Q)\hookrightarrow H^{1/2}(\G,\C^2)\hookrightarrow L^{p}(\G,\C^{2}), \qquad\forall p\geq 2.
\]
On the other hand, according to \eqref{eq-sp_D} one can decompose the form domain as the orthogonal sum of the positive and negative \emph{spectral subspaces} for the operator $\Dg$, i.e.
\[
 \d(\Q)=\d(\Q)^{+}\oplus \d(\Q)^{-}.
\]
As a consequence, if one denotes $\psi^+:=P^+\psi$ and $\psi^{-}:=P^-\psi$ and recalls that it is possible to define a norm for $\d(\Q)$ as
\[
\|\psi\|_{\Q}:=\|\sqrt{|\Dg|}\psi\|_{2,\G},\qquad \forall\psi\in \d(\Q)
\]
(with $\sqrt{|\Dg|}\psi$ given by Borel functional calculus for $\Dg$), then there results that 
\begin{equation}
 \label{eq-qexpl}
 \Q(\psi)=\frac{1}{2}\left(\|\psi^+\|_{\Q}-\|\psi^-\|_{\Q}\right).
\end{equation}

Clearly, \eqref{eq-domex} and \eqref{eq-qexpl} are not explicit forms for the domain and the action of the quadratic form associated with $\Dg$. Nevertheless, they present the suitable detalis in order to manage the computations required in the proofs of Theorems \ref{thm-bound} and \ref{thm-limit} below.


\subsection{Bound states}

The study of bound states of \eqref{eq-NLDconc} presents some relevant differences with respect to the NLS case \eqref{eq-NLSconc} (a difference which also arises in the extended case).

The main point here is the unboundedness from below of the spectrum of $\Dg$ which makes the associated quadratic form \emph{strongly indefinite}, even fixing the $L^2$-norm. As a consequence, the natural energy functional associated with \eqref{eq-NLDconc}, i.e.
\[
 \Q(\psi)-\frac1p\int_\K|\psi|^{p-2}\psi,\qquad p\geq2,
\]
is unbounded from below even under the mass constraint $\|\psi\|_{2,\G}^2=\mu$. Hence, no minimization can be performed and also the search for constrained critical points presents several technical difficulties.

Therefore, the most promising strategy seems to be that of considering the \emph{action functional} associated with \eqref{eq-NLDconc} without any constraint. Then for a fixed $\omega\in\R$ one looks for critical points of the functional
\[
 \L(\psi,\K,p)=\Q(\psi)-\f{\omega}{2}\int_{\G}|\psi|^{2}-\frac{1}{p}\int_{\K}|\psi|^{p}\dx,
\]
on the form-domain $\d(\Q)$.

To this aim, one needs to adapt to the metric graphs setting some well known techniques of Critical Point Theory. In particular, a suitable version of the so-called \emph{linking} technique for even functional is exploited. We refer the reader to \cite{R-CBMS,S}, for a general presentation of those methods in an abstract setting, together with several applications, and to \cite{ES-CMP} which deals with nonlinear Dirac equations.

\medskip
Critical levels for $\L$ are defined as the following min-max levels
\[
\alpha_{N}:=\inf_{X\in\F_{N}}\sup_{\psi\in X}\L(\psi,\K,p),
\]
with
\[
\F_{N}:=\left\{X\in H^1(\G,\C^2)\backslash\{0\}:X\text{ closed and symmetric, s.t. }\gamma[h_{t}(S^{+}_{r})\cap X]\geq N,\,\forall t\geq0 \right\},
\]
where $h_t$ is the usual \emph{pseudo-gradient flow} associated with $\L(\cdot,\K,p)$, $S^{+}_{r}$ is the sphere of radius $r$ of $\d(\Q)^+$ and $\gamma$ denotes again the Krasnosel'skii genus. Indeed, one can check that, if $\F_N\neq\emptyset$, then there exists a Palais-Smale sequence $\left( \psi_{n}\right)\subset \d(\Q)$ at level $\alpha_{N}$, i.e. 
\[
 \begin{cases}
 \displaystyle \L(\psi_{n},\K,p)\to \alpha_{N}  \\[.2cm]
 \displaystyle d\L(\psi_{n},\K,p)\xrightarrow{\d(\Q)'} 0,
 \end{cases}
\]
and in addition, there results
\begin{gather*}
\alpha_{N_{1}}\leq\alpha_{N_{2}}, \qquad \forall N_{1}<N_{2},\\[.2cm]
\rho\leq\alpha_{N}\leq \widetilde{\rho}<+\infty,\qquad\forall N\in\mathbb{N},
\end{gather*}
for fixed $\rho,\,\widetilde{\rho}>0$.

\begin{figure}[t]
\centering
{\includegraphics[width=.55\columnwidth]{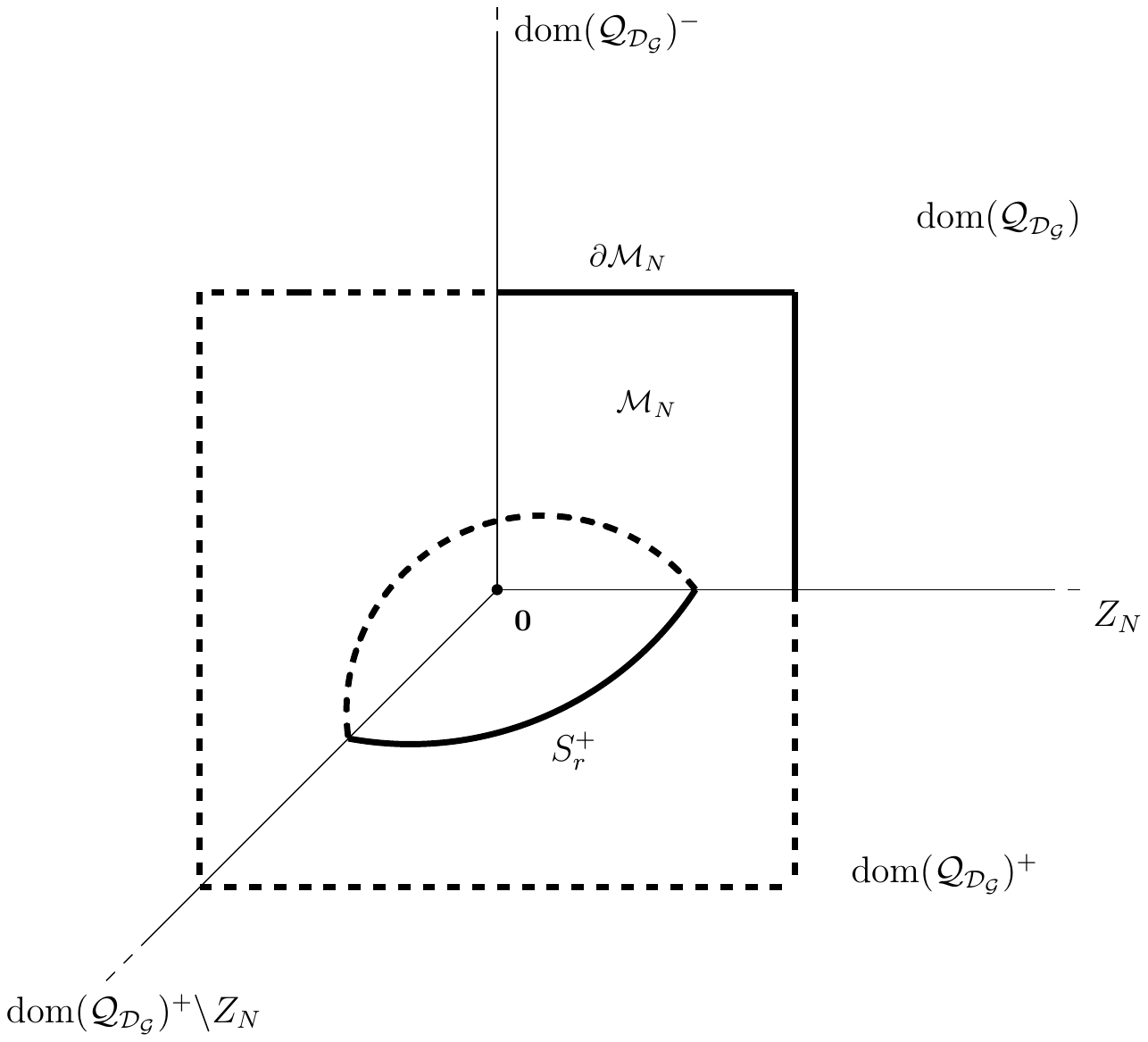}}
\caption{A graphic insight on linking geometry.}
\label{fig-linking}
\end{figure}

As a consequence, in order to prove that $\F_N\neq\emptyset$, it suffices to show that the functional possesses a so-called \emph{linking geometry}; namely, that for every $N\in N$ there exist $R=R(N,p)>0$ and an $N$-dimensional space $Z_{N}\subset \d(\Q)^{+}$ such that
\[
\L(\psi,\K,p)\leq 0,\qquad \forall\psi\in\partial\M_N,
\]
where
\[ 
\partial \M_{N}=\left\{\psi\in\M_N : \|\psi^{-}\|=R \text{ or }\|\psi^{+}\|=R \right\}
\]
and
\[
\M_{N}:=\left\{\psi\in \d(\Q) : \|\psi^{-}\|\leq R \text{ and } \psi^{+}\in Z_{N} \text{ with }\|\psi^{+}\|\leq R \right\},
\]
and that there exist $r,\rho>0$ such that 
\[
\inf_{S^{+}_{r}}\L(\cdot,\K,p)\geq\rho>0
\]
(for a graphic intuition see, e.g., Figure \ref{fig-linking}).

Finally, checking the validity of the Palais-Smale condition at positive levels for the action functional it is possible to prove the following

\begin{thm}[{\cite[Theorem 2.11]{BCT-p}}]
 \label{thm-bound}
 Let $\G$ satisfy {\bf (H1)-(H3)}, $p>2$ and $m,c>0$. Then, for every $\omega\in\R\backslash\sigma(\Dg)=(-mc^2,mc^2)$ there exists infinitely many (distinct) pairs of bound states of frequency $\omega$ of the NLDE, at mass $m$ and relativistic parameter $c$.
\end{thm}


\section{Nonrelativistic limit}
\label{sec-limit}

The NLDE and the NLSE equation are clearly closely related as, physically, the latter should correspond to the nonrelativistic limit of the former. Heuristically, one expect to recover the NLSE from NLDE as the relativistic parameter tends to infinity (namely, ``$c\to\infty$''), that is, neglecting relativistic effects. It is, then, particularly interesting to rigorously establish such a connection. 

This has been first done in \cite{ES-ANHIPA} for the Dirac-Fock equations proving the convergence of bound states to those of the (nonrelativistic) Hartree-Fock model.

\medskip
More in detail, given a sequence $c_n\to\infty$, one is interested in the limit behavior of a sequence of bound states with frequencies $0<\omega_n<mc^2_n$. Precisely, one has to choose frequencies such that
\[
 \omega_n-mc^2_n\longrightarrow\frac{\lambda}{m}<0.
\]
This represents a kind of renormalization, which corresponds to the fact that one has to subtract the \emph{rest energy} of the particles, as this property is a peculiarity of relativistic theories which is absent in the nonrelativistic setting.

The strategy adopted in \cite{BCT-p} is an adaptation of the one developed in \cite{ES-ANHIPA}. Namely, first one has to establish $H^1$ uniform bounds for the sequence of the bound states and then, suitably manipulating \eqref{eq-NLDconc}, one has to prove that the lower component goes to zero while the sequence of upper components is a Palais-Smale sequence for the action functional associated to \eqref{eq-NLSconc}, namely for 
\[
 \mathcal{J}(u):=\f{1}{2}\int_\G|u'|^2\dx-\frac{2m}{p}\int_\K|u|^p\dx-\f{\lambda}{2}\int_\G|u|^2\dx.
\]
at a fixed $\lambda<0$. Therefore, up to a compactness argument, the existence of a limit satisfying the NLSE is proved.

However, we would like to remark that some relevant differences are present here with respect to \cite{ES-ANHIPA}, that call for some relevant modifications of the proofs. In particular, in \cite{ES-ANHIPA} bound states were obtained as constrained critical points of the energy. As a consequence uniform boundedness of the $L^2$-norms and the non-triviality limit are easily obtained. On the contrary, in the case considered in \cite{BCT-p}, those properties are not a priori guaranteed and represent two main points of the proofs.

However, for the sake of completeness, we mention that in \cite{ES-ANHIPA} one of the most delicate parts the proof is to estimate the sequence of the Lagrange multipliers corresponding to the $L^2$-constraint, proving that they are in the spectral gap of the operator. This, actually, is one of the reasons for which in \cite{BCT-p} we looked for critical points of the action and not for the constrained critical points of the energy, in order to avoid such additional technical difficulties.

In view of all the above remarks, it is possible to state the following

\begin{thm}[{\cite[Theorem 2.12]{BCT-p}}]
 \label{thm-limit}
 Let $\G$ satisfy {\bf (H1)-(H3)}, $p\in(2,6)$, $m>0$ and $\lambda<0$. Let also $(c_{n}), (\omega_{n})$ be two real sequences such that
 \begin{gather*}
 0<c_{n},\omega_n\rightarrow+\infty,\\[.2cm]
 \omega_{n}<mc^{2}_{n},\\[.2cm]
 \omega_n-mc_n^2\longrightarrow\f{\lambda}{m}.
\end{gather*}
If $\{\psi_n=(\phi_n,\chi_n)^T\}$ is a bound state of frequency $\omega_n$ of the NLDE at relativistic parameter $c_n$, then, up to subsequences, there holds
 \[
  \phi_n\to u\qquad\text{and}\qquad\chi_n\to0\qquad\text{in}\quad H^1(\G).
 \]
 where $u$ is a bound state of frequency $\lambda$ of the NLSE.
\end{thm}

Actually, the bound states of NLDE do not converge exactly to the bound states of the NLSE depicted in \eqref{eq-NLDconc}. In fact, $u$ is a solution of 
\[
 -\lap u-2m\chi_{\K}|u|^{p-2}\,u=\lambda u.
\]
The coefficient $2m$ is consistent with the fact that in the nonrelativistic limit the kinetic (free) part of the hamiltonian of a particle is given by $$\frac{-\Delta}{2m}\,.$$

\begin{rem}
 It is also worth observing that Theorem \ref{thm-limit}, in contrast to Theorem \ref{thm-bound}, holds only for a fixed range of exponents, the $L^2$-subcritical case $p\in(2,6)$.
\end{rem}



\begin{thebibliography}{99}

\bibitem{ACFN-RMP}
R. Adami, C. Cacciapuoti, D. Finco, D. Noja,
Fast solitons on star graphs,
\emph{Rev. Math. Phys.} {\bf 23} (2011), no. 4, 409--451.

\bibitem{ACFN-JPA}
R. Adami, C. Cacciapuoti, D. Finco, D. Noja,
On the structure of critical energy levels for the cubic focusing NLS on star graphs,
\emph{J. Phys. A} {\bf 45} (2012), no. 19, article number 192001, 7pp.

\bibitem{ACFN-JDE}
R. Adami, C. Cacciapuoti, D. Finco, D. Noja,
Variational properties and orbital stability of standing waves for NLS equation on a star graph,
\emph{J. Differential Equations} {\bf 257} (2014), no. 10, 3738--3777.

\bibitem{ACFN-ANHIPC}
R. Adami, C. Cacciapuoti, D. Finco, D. Noja,
Constrained energy minimization and orbital stability for the NLS equation on a star graph,
\emph{Ann. Inst. H. Poincaré Anal. Non Lin\'eaire} {\bf 31} (2014), no. 6, 1289--1310.

\bibitem{AM}
A. Ambrosetti A., A. Malchiodi,
\emph{Nonlinear analysis and semilinear elliptic problems},
Cambridge Studies in Advanced Mathematics 104, Cambridge University Press, Cambridge, 2007.

\bibitem{ADST-APDE}
R. Adami, S.Dovetta, E. Serra, P. Tilli, 
Dimensional crossover with a continuum of critical exponents for NLS on doubly periodic metric graphs, 	arXiv:1805.02521 [math.AP] (2018). Accepted by \emph{Anal. PDE}.

\bibitem{AST-CVPDE}
R. Adami, E. Serra, P. Tilli,
NLS ground states on graphs,
\emph{Calc. Var. Partial Differential Equations} {\bf 54} (2015), no. 1, 743--761.

\bibitem{AST-JFA}
R. Adami, E. Serra, P. Tilli,
Threshold phenomena and existence results for NLS ground states on metric graphs,
\emph{J. Funct. Anal.} {\bf 271} (2016), no. 1, 201--223.  

\bibitem{AST-CMP}
R. Adami, E. Serra, P. Tilli, 
Negative energy ground states for the L2-critical NLSE on metric graphs,
\emph{Comm. Math. Phys.} {\bf 352} (2017), no. 1, 387--406.

\bibitem{AST-PARMA}
R. Adami, E. Serra, P. Tilli,
Nonlinear dynamics on branched structures and networks.,
\emph{Riv. Math. Univ. Parma (N.S.)} {\bf 8} (2017), no. 1, 109--159.

\bibitem{AST-CVPDE2}
R. Adami, E. Serra, P. Tilli,
Multiple positive bound states for the subcritical NLS equation on metric graphs,
\emph{Calc. Var. Partial Differential Equations} {\bf 58} (2019), no. 1, article number 5, 16pp.

\bibitem{AS-JMP}
J. Arbunich, C. Sparber,
Rigorous derivation of nonlinear Dirac equations for wave propagation in honeycomb structures,
\emph{J. Math. Phys.} {\bf 59} (2018), no. 1, article number 011509, 19pp.

\bibitem{BL-ARMA}
H. Berestycki, P.-L. Lions,
Nonlinear scalar field equations II. Existence of infinitely many solutions,
\emph{Arch. Rational Mech. Anal.} {\bf 82} (1983), no. 4, 347--375.

\bibitem{BK}
G. Berkolaiko, P. Kuchment,
\emph{Introduction to quantum graphs},
Mathematical Surveys and Monographs 186, American Mathematical Society, Providence, RI, 2013.

\bibitem{B-JDE}
W. Borrelli,
Stationary solutions for the 2D critical Dirac equation with Kerr nonlinearity,
\emph{J. Differential Equations} {\bf 263} (2017), no. 11, 7941--7964.

\bibitem{B-JMP}
W. Borrelli,
Multiple solutions for a self-consistent Dirac equation in two dimensions,
\emph{J. Math. Phys.} {\bf59} (2018), no. 4, article number 041503, 13pp.

\bibitem{B-CVPDE}
W. Borrelli,
Weakly localized states for nonlinear Dirac equations,
\emph{Calc. Var. Partial Differential Equations} {\bf 57} (2018), no. 6, article number 155, 21pp.

\bibitem{BCT-p}
W. Borrelli, R. Carlone, L. Tentarelli,
Nonlinear Dirac equation on graphs with localized nonlinearities: bound states and nonrelativistic limit,
arXiv:1807.06937 [math.AP] (2018).

\bibitem{BT-JMP}
W. Bulla, T. Trenkler,
The free Dirac operator on compact and noncompact graphs, 
\emph{J. Math. Phys.} {\bf 31} (1990), no. 5, 1157--1163. 

\bibitem{CDS-MJM}
C. Cacciapuoti, S. Dovetta, E. Serra,
Variational and stability properties of constant solutions to the NLS equation on compact metric graphs,
\emph{Milan J. Math.} {\bf 86} (2018), no. 2, 305--327.

\bibitem{CFN-PRE}
C. Cacciapuoti, D. Finco, D. Noja,
Topology-induced bifurcations for the nonlinear Schr\"odinger equation on the tadpole graph,
\emph{Phys. Rev. E (3)} {\bf 91} (2015), no. 1, article number 013206, 8pp.

\bibitem{CFN-Non}
C. Cacciapuoti, D. Finco, D. Noja,
Ground state and orbital stability for the NLS equation on a general starlike graph with potentials,
\emph{Nonlinearity} {\bf 30} (2017), no. 8, 3271--3303.

\bibitem{C}
T. Cazenave,
\emph{Semilinear Schr\"odinger equations},
Courant Lecture Notes in Mathematics 10, American Mathematical Society, Providence, RI, 2003.

\bibitem{D-JDE}
S. Dovetta,
Existence of infinitely many stationary solutions of the $L^2$-subcritical and critical NLSE on compact metric graphs,
\emph{J. Differential Equations} {\bf 264} (2018), no. 7, 4806--4821.

\bibitem{DD-p}
S. Dovetta,
Mass-constrained ground states of the stationary NLSE on periodic metric graphs,
arXiv:1811.06798 [math.AP] (2018).

\bibitem{DT-OTAA}
S. Dovetta, L. Tentarelli,
Ground states of the $L^2$-critical NLS equation with localized nonlinearity on a tadpole graph,
arXiv:1803.09246 [math.AP] (2018). Accepted by \emph{Oper. Theory Adv. Appl.}

\bibitem{DT-p}
S. Dovetta, L. Tentarelli,
$L^2$-critical NLS on noncompact metric graphs with localized nonlinearity: topological and metric features, arXiv:1811.02387 [math.AP] (2018).

\bibitem{D-p}
A. Duca,
Global exact controllability of the bilinear Schr\"odinger potential type models on quantum graphs,
arXiv:1710.06022 [math.OC] (2017).

\bibitem{ES-CMP}
M.J. Esteban, E. S\'er\'e,
Stationary states of the nonlinear Dirac equation: a variational approach,
\emph{Comm. Math. Phys.} {\bf 171} (1995), no. 2, 323--350.

\bibitem{ES-ANHIPA}
M.J. Esteban, E. S\'er\'e,
Nonrelativistic limit of the Dirac-Fock equations,
\emph{Ann. Henri Poincar\'e} {\bf 2} (2001), no. 5, 941--961. 

\bibitem{FW-CMP}
C.L. Fefferman, M.I. Weinstein,
Wave Packets in Honeycomb Structures and Two-Dimensional Dirac Equations,
\emph{Comm. Math. Phys.} {\bf 326} (2014), no. 1, 251--286. 

\bibitem{GPS-NODEA}
S. Gilg, D. Pelinovsky, G. Schneider,
Validity of the NLS approximation for periodic quantum graphs, 
\emph{NoDEA Nonlinear Differential Equations Appl.} {\bf 23} (2016), no. 6, article number 63, 30pp.

\bibitem{GSD-PRA}
S. Gnutzmann, U. Smilansky, S. Derevyanko,
Stationary scattering from a nonlinear network,
\emph{Phys. Rev. A} {\bf 83} (2011), no. 3, article number 033831, 6pp.

\bibitem{GW-PRE}
S. Gnutzmann, D. Waltner,
Stationary waves on nonlinear quantum graphs: general framework and canonical perturbation theory,
\emph{Phys. Rev. E} {\bf 93} (2016), no. 3, article number 032204, 19pp.

\bibitem{HC-PD}
L.H. Haddad, L.D. Carr,
The nonlinear Dirac equation in Bose-Einstein condensates: foundation and symmetries,
\emph{Phys. D} {\bf 238} (2009), no. 15, 1413--1421.

\bibitem{KP-JDE}
A. Kairzhan, D.E. Pelinovsky,
Nonlinear instability of half-solitons on star graphs,
\emph{J. Differential Equations} {\bf 264} (2018), no. 12, 7357--7383.

\bibitem{K}
M. A. Krasnosel'skii,
\emph{Topological methods in the theory of nonlinear integral equations},
A Pergamon Press Book, The Macmillan Co., New York, 1964.

\bibitem{MP-AMRX}
J.L. Marzuola, D.E. Pelinovsky,
Ground state on the dumbbell graph,
\emph{Appl. Math. Res. Express. AMRX} (2016), no. 1, 98--145. 

\bibitem{MNS-APDE}
D. Mugnolo, D. Noja, C. Seifert,
Airy-type evolution equations on star graphs,
\emph{Anal. PDE} {\bf 11} (2018), no. 7, 1625--1652.

\bibitem{N-RSTA}
D. Noja, 
Nonlinear Schr\"odinger equation on graphs: recent results and open problems,
\emph{Philos. Trans. R. Soc. Lond. Ser. A Math. Phys. Eng. Sci.} {\bf 372} (2014), no. 2007, article number 20130002, 20pp.

\bibitem{NPS-Non}
D. Noja, D. Pelinovsky, G. Shaikhova, 
Bifurcations and stability of standing waves in the nonlinear Schr\"odinger equation on the tadpole graph,
\emph{Nonlinearity} {\bf 28} (2015), no. 7, 2343--2378.

\bibitem{NRS-JDE}
D. Noja, S. Rolando, S. Secchi,
Standing waves for the NLS on the double-bridge graph and a rational-irrational dichotomy,
{\em J. Differential Equations} {\bf 266} (2019), no. 1, 147--178.

\bibitem{PS-ANIHPA}
D. Pelinovsky, G. Schneider,
Bifurcations of standing localized waves on periodic graphs, 
\emph{Ann. Henri Poincar\'e} {\bf 18} (2017), no. 4, 1185--1211.

\bibitem{P}
O. Post,
Equilateral quantum graphs and boundary triples,
\emph{Analysis on graphs and its applications} 469--490, 
Proc. Sympos. Pure Math. 77, AMS, Providence, RI, 2008. 

\bibitem{R-CBMS}
P.H. Rabinowitz,
\emph{Minimax methods in critical point theory with applications to differential equations},
CBMS Regional Conference Series in Mathematics 65, AMS, Providence, RI, 1986.

\bibitem{SBMK-JPA}
K.K. Sabirov, D.B. Babajanov, D.U. Matrasulov, P.G. Kevrekidis,
Dynamics of Dirac solitons in networks,
\emph{J. Phys. A} {\bf 51} (2018), no. 43, article number 435203, 13pp.

\bibitem{ST-JDE}
E. Serra, L. Tentarelli, 
Bound states of the NLS equation on metric graphs with localized nonlinearities,
\emph{J. Differential Equations} {\bf 260} (2016), no. 7, 5627--5644.

\bibitem{ST-NA}
E. Serra, L. Tentarelli,
On the lack of bound states for certain NLS equations on metric graphs,
\emph{Nonlinear Anal.} {\bf 145} (2016), 68--82.

\bibitem{S}
M. Struwe,
\emph{Variational methods. Applications to nonlinear partial differential equations and Hamiltonian systems},
Fourth edition, Results in Mathematics and Related Areas, 3rd Series, A Series of Modern Surveys in Mathematics, 34, Springer-Verlag, Berlin, 2008.

\bibitem{T-JMAA}
L. Tentarelli, NLS ground states on metric graphs with localized nonlinearities,
\emph{J. Math. Anal. Appl.} {\bf 433} (2016), no. 1, 291--304.

\bibitem{TLB-AP}
T.X. Tran, S. Longhi, F. Biancalana,
Optical analogue of relativistic Dirac solitons in binary waveguide arrays,
\emph{Ann. Physics} {\bf 340} (2014), 179--187.

\end{thebibliography}
\end{document}